\documentclass{article}

\usepackage{graphicx}
\usepackage{amsmath}
\usepackage{amssymb}
\usepackage{url}

\newtheorem{theorem}{Theorem}[section]
\newtheorem{proposition}[theorem]{Proposition}
\newtheorem{remark}[theorem]{Remark}

\newcommand{\bm}[1]{\mathbf{#1}}
\newcommand{\A}{\mathcal{A}}
\newcommand{\B}{\mathcal{B}}
\newcommand{\e}{{\mathrm{e}}}
\newcommand{\ii}{{\mathbf{i}}}
\newcommand{\N}{{\mathbb{N}}}
\newcommand{\R}{{\mathbb{R}}}
\newcommand{\C}{{\mathbb{C}}}
\newcommand{\cD}{{\mathcal{D}}}
\newcommand{\cS}{{\mathcal{S}}}
\newcommand{\el}{\ensuremath{\textnormal{{\boldmath$\ell$}}}}
\newcommand{\sigmab}{\boldsymbol\sigma}
\newcommand{\alphab}{\boldsymbol\alpha}
\newcommand{\omegab}{\boldsymbol\omega}
\newcommand{\phib}{\boldsymbol\phi}
\newcommand{\psib}{\boldsymbol\psi}
\newcommand{\flops}{\emph{flops}}
\DeclareMathOperator{\diag}{diag}
\DeclareMathOperator{\cond}{cond}
\newcommand{\ttt}[1]{\texttt{#1}}

\usepackage{ifthen}
\usepackage{float}
\floatstyle{ruled}
\newfloat{algorithm}{htb}{alg}
\floatname{algorithm}{Algorithm}
\newcounter{algo@row}
\newcounter{algo@rowindent}
\newcommand{\algofont}[1]{\textbf{#1}}
\newcommand{\algonumbersize}[1]{\scriptsize{#1}}
\newcommand{\algopreitem}[1][\arabic{algo@row}]{\texttt{\algonumbersize{#1}}}
\newcommand{\algoitemskip}{\hspace{\value{algo@rowindent}cc}}
\newcommand{\algosetline}[1]{\setcounter{algo@row}{#1}\addtocounter{algo@row}{-1}}
\newcommand{\algoskipindent}[1]{\addtocounter{algo@rowindent}{#1}}
\newenvironment{algo}{\vskip.3em\small%
  \begin{list}{\algopreitem\texttt{\algonumbersize{:}}}{%
      \usecounter{algo@row}%
      \setcounter{algo@rowindent}{0}%
      \setlength{\itemindent}{2em}%
      \setlength{\labelwidth}{2em}
      \setlength{\parsep}{0cm}%
    }%
}{
  \end{list}\vskip-.5em
}
\newcommand{\algonewnestedopen}[2]{
  \newcommand{#1}[1][]{%
    \ifthenelse{\equal{##1}{}}{\item}{\item[{\algopreitem[##1]}]}
    \algoitemskip\algofont{#2}%
    \addtocounter{algo@rowindent}{1}%
    \ignorespaces
  }
}
\newcommand{\algonewnestedaux}[2]{
  \newcommand{#1}[1][]{
    \addtocounter{algo@rowindent}{-1}
    \ifthenelse{\equal{##1}{}}{\item}{\item[{\algopreitem[##1]}]}
    \algoitemskip\algofont{#2}%
    \addtocounter{algo@rowindent}{+1}%
    \ignorespaces
  }
}
\newcommand{\algonewnestedclose}[2]{
  \newcommand{#1}[1][]{
    \addtocounter{algo@rowindent}{-1}
    \ifthenelse{\equal{##1}{}}{\item}{\item[{\algopreitem[##1]}]}
    \algoitemskip\algofont{#2}%
    \ignorespaces
  }
}
\newcommand{\algonewcommand}[2]{
  \newcommand{#1}[1][default]{
    \ifthenelse{\equal{##1}{default}}{\item}{\item[{\algopreitem[##1]}]}%
    \algoitemskip\algofont{#2}%
    \ignorespaces
  }%
}
\newcommand{\algonewkeyword}[2]{\newcommand{#1}{\algofont{#2}}}
\algonewcommand{\STATE}{\ignorespaces}
\algonewcommand{\REQUIRE}{Require: }
\algonewcommand{\COMPUTE}{Compute: }
\algonewcommand{\OUTPUT}{Output: }
\newcommand{\COMMENT}{\hfill // }
\algonewnestedopen{\IF}{if }
\algonewnestedaux{\ELSEIF}{else if }
\algonewnestedclose{\ENDIF}{end if }
\algonewnestedopen{\FOR}{for }
\algonewnestedclose{\ENDFOR}{end for }
\algonewkeyword{\To}{to }%
\algonewkeyword{\If}{if }%
\algonewkeyword{\Then}{then }%

\title{A fast solver for linear systems \\ with displacement structure}

\author{Antonio Aric\`o\thanks{Dipartimento di Matematica e Informatica,
Universit\`a di Cagliari, viale Merello 92, 09123 Cagliari, Italy. E-mail:
\texttt{\{arico,rodriguez\}@unica.it}. This work was partially supported by
MIUR, under the PRIN grant no.~20083KLJEZ-003, and by INdAM-GNCS.} 
\and Giuseppe Rodriguez\footnotemark[1]}

\date{}

\begin{document}

\maketitle

\begin{abstract} 
We describe a fast solver for linear systems with reconstructable Cauchy-like
structure, which requires $O(rn^2)$ floating point operations and $O(rn)$
memory locations, where $n$ is the size of the matrix and $r$ its
\emph{displacement rank}.
The solver is based on the application of the generalized Schur algorithm 
to a suitable augmented matrix, under some assumptions on the knots of the
Cauchy-like matrix.
It includes various pivoting strategies, already discussed in the literature,
and a new algorithm, which only requires reconstructability.
We have developed a software package, written in Matlab and C-MEX, which
provides a robust implementation of the above method.
Our package also includes solvers for Toeplitz(+Hankel)-like and
Vandermonde-like linear systems, as these structures can be reduced to
Cauchy-like by fast and stable transforms.
Numerical experiments demonstrate the effectiveness of the software.
\end{abstract}

%
%
\section{Introduction}\label{sec:intro}
%
%
The idea to connect the solution of some linear problems to suitable
augmented matrices has often been used in linear algebra; see, 
e.g.,~\cite{bgs70,bjo92,sie65} for an application to full rank least
squares problems. More recently, it was proposed~\cite{kc94} to solve
various computational problems involving Toeplitz matrices by 
evaluating a Schur complement of a suitable augmented matrix.

The reason for this approach lies in the fact that Toeplitz matrices,
as well as other classes of structured matrices, can be expressed as
the solution of a \emph{displacement equation}~\cite{fmkl79,kkm79b},
and that this property is inherited by their Schur complements of any
order. This fact allows one to store a structured matrix of dimension $n$
in $O(n)$ memory locations and, which is most important, 
to apply a fast implementation of
the Gauss triangularization method (the \emph{generalized Schur
algorithm}) operating directly on the displacement information
associated to the matrix~\cite{gko95,hei95a,hb98}; see also~\cite{ks95,ks99}
for a review.

In this paper we use the mathematical tools above mentioned to devise
a fast algorithm to solve the linear system
\begin{equation}\label{eq:Cx=b}
  C\bm{x}=\bm{b},
\end{equation}
where $\bm{b}\in\C^{n}$ and $C\in\C^{n\times n}$ is a Cauchy-like matrix;
this means that
\begin{equation}\label{eq:Cl}
  C_{ij} = 
    \dfrac{\phib_i^*\,\psib_j}{t_i-s_j},
\end{equation}
where $\phib_i,\psib_j\in\C^r$, $t_i,s_j\in\C$, $i,j=1,\dots,n$, and
the asterisk denotes the conjugate transpose.

Given a matrix
\begin{equation*}
  M = \begin{bmatrix}
    A & B \\
    C & D
  \end{bmatrix},
\end{equation*}
with $A\in\C^{n\times n}$ invertible, we denote its \emph{Schur
  complement of order $n$} by
\begin{equation*}
  \cS_n(M)=D-CA^{-1}B.
\end{equation*}
We associate to the system~\eqref{eq:Cx=b} the augmented matrix
\begin{equation}\label{eq:ACb}
  \A_{C,\bm{b}} = \begin{bmatrix}C&\bm{b}\\-I_n&O\end{bmatrix},
\end{equation}
where $I_n$ is the identity matrix of size $n$ and $O$ is a null matrix.
Note that $\A_{C,\bm{b}}$ is Cauchy-like too; see
Section~\ref{sec:displ}. We will apply the generalized Schur
algorithm (GSA) to $\A_{C,\bm{b}}$, to compute the solution of the
system~\eqref{eq:Cx=b} as the Schur complement
\[
\cS_n(\A_{C,\bm{b}})=C^{-1}\bm{b}.
\]
A similar approach was used in~\cite{rod06} to solve structured least
square problems.

It is remarkable that other classes of structured linear systems can
be reduced to Cauchy-like systems by fast transforms (see~\cite{gko95,hb97}),
for example the systems whose matrix is either Toeplitz,
Hankel, or Vandermonde. So, this approach is quite general.

As observed in~\cite{gko95,hei95a}, the advantage of applying the GSA
to a Cauchy-like matrix is that such a structure is invariant under
rows/columns exchanges, so it is easy to embed a pivoting strategy in
the algorithm.
The problem of choosing a pivoting technique for Cauchy matrices is addressed
in~\cite{cr97,gu98,sb95}.

Our approach requires $O(n^2)$ floating
point operations and $O(n)$ memory locations. On the contrary, the
approach adopted in~\cite{gko95} consists of explicitly
computing (and storing) the LU factorization of $C$, hence requiring
$O(n^2)$ locations.

We developed a software package, called \ttt{drsolve}, which includes a
Cauchy-like solver, some interface routines for other structured linear
systems, and various conversion and auxiliary routines.
The package is written in Matlab~\cite{matlab79} for the most part. The
Cauchy-like solver has also been implemented in C language, with
extensive use of the BLAS library~\cite{blas02}, and it has been linked to
Matlab via the MEX (\emph{Matlab executable}) interface library. The software
includes a device to detect numerical singularity or ill-conditioning
of the coefficient matrix.

Among the other software for structured linear systems publicly available, we
would like to mention the Toeplitz Package~\cite{toepack83}, written in a
Russian-American collaboration, and \ttt{toms729}~\cite{hc92}, based on a
look-ahead Levinson algorithm, which we will use in our numerical experiments.
Various computer programs for structured problems have been developed by the
MaSe-team~\cite{maseteam}, coordinated by Marc Van Barel, and by the members of
our research group~\cite{bezout}, coordinated by Dario Bini.

This paper is organized as follows: in Section~\ref{sec:displ} we recall the
displacement equations of some classes of structured matrices, and the rules
to convert them to Cauchy-like form.
In Sections~\ref{sec:GKO4A} and~\ref{sec:pivprec}, we describe the algorithms
for solving a Cauchy-like linear system, and the pivoting techniques,
respectively, that we have implemented.
In Section~\ref{sec:impl} we give some details on the software package, and in
Section~\ref{sec:numres} we discuss the results of a widespread numerical
experimentation.

%
%
\section{Displacement structure}\label{sec:displ}
%
%
A matrix $A\in\C^{n\times n}$ is said to satisfy a Sylvester
displacement equation if
\begin{equation}\label{eq:EAAF=}
  EA-AF=GH^*,
\end{equation}
where $G,H\in\C^{n\times r}$ are full rank matrices, called the
generators, $E,F\in\C^{n\times n}$ are the displacement matrices, and
$r$ is the displacement rank~\cite{fmkl79,hr84,kkm79b}. This
representation is particularly relevant when $r$ is
significantly smaller than $n$.

A Cauchy-like matrix~\eqref{eq:Cl} satisfies
equation~\eqref{eq:EAAF=} with $E=D_\bm{t}=\diag(t_1,\dots,t_n)$,
$F=D_\bm{s}=\diag(s_1,\dots,s_n)$, and
\begin{equation}\label{eq:GHC}
  G=\begin{bmatrix}\phib_1&\cdots&\phib_n\end{bmatrix}^*,
  \qquad
  H=\begin{bmatrix}\psib_1&\cdots&\psib_n\end{bmatrix}^*.
\end{equation}
Any matrix $A$ satisfying~\eqref{eq:EAAF=} can be converted to a
Cauchy-like one if both $E$ and $F$ are diagonalizable. In fact, if $E=U
D_\bm{t} U^{-1}$ and $F = V D_\bm{s} V^{-1}$, then~\eqref{eq:EAAF=}
becomes
\begin{equation*}
  D_\bm{t} C - C D_\bm{s} = \widetilde{G}\,\widetilde{H}^*,
\end{equation*}
where $C = U^{-1} A V$, $\widetilde{G} = U^{-1} G$, and $\widetilde{H}
= V^* H$. The same transformation converts a linear system
$A\bm{x}=\bm{b}$ to $C\widetilde{\bm{x}}=\widetilde{\bm{b}}$, where
$\widetilde{\bm{b}}=U^{-1}\bm{b}$ and $\bm{x}=V\widetilde{\bm{x}}$.
This transformation is numerically effective for large matrices
if $U$ and $V$ are unitary and the computation is fast.

\begin{table}
\small
\begin{equation*}
  \begin{array}{l|cc|cccc}
    \multicolumn{1}{c|}{\textnormal{structure}} & E & F & U & V & D_\bm{t} &
    D_\bm{s} \\
    \hline
    \textnormal{Cauchy-like} & D_\bm{t} & D_\bm{s} \\
    \textnormal{Toeplitz-like} & Z_1 & Z_{-1} & F_1 & F_{-1} & D_1 & D_{-1} \\
    \textnormal{Toeplitz+Hankel-like} & Y_0 & Y_1 & \mathcal{S} & \mathcal{C} &
	    D_\mathcal{S} & D_\mathcal{C} \\
    \textnormal{Vandermonde-like} & D_\bm{w} & Z_\phi^{*} & I_n & F_\phi & 
	    D_\bm{w} & D_\phi^* \\
    \hline
  \end{array}
\end{equation*}
\caption{Displacement matrices for some structured classes.}
\label{tab:struct}
\end{table}
Some well-known classes of structured matrices which fall within this
framework are reported in Table~\ref{tab:struct}; see~\cite{gko95,hei95a,hb97}. 
The corresponding displacement matrices are defined as follows:
\begin{equation*}
  D_\bm{v} = \diag(v_1,\dots,v_n),
  \quad
  Z_\phi =
  \begin{bmatrix}
    & \!\phi\\
    I_{n-1}
  \end{bmatrix},
  \quad
  Y_\delta = 
  \begin{bmatrix}
    \,\,\delta\,\, & 1      & 0      & \hdots & 0             \\[-1mm]
    1              & 0      & 1      & \ddots & \vdots        \\[-1mm]
    0              & 1      & \ddots & \ddots & 0             \\[-1mm]
    \vdots         & \ddots & \ddots & 0      & 1             \\
    0              & \hdots & 0      & 1      & \,\,\delta\,\,
  \end{bmatrix},
\end{equation*}
where $|\phi|=1$ and $\delta=0,1$. 
Their spectral factorization is
analytically known, as
\begin{equation*}
  Z_\phi = F_\phi D_\phi F_\phi^*,
  \qquad
  Y_0 = \mathcal{S} D_\mathcal{S} \mathcal{S},
  \qquad 
  Y_1 = \mathcal{C} D_\mathcal{C} \mathcal{C}^T,
\end{equation*}
where, setting $\omega=\e^{\frac{2\pi\ii}n}$, $q_1=2^{-1/2}$, $q_\ell=1$,
$\ell=2,\dots,n$, and denoting by $\phi^{1/n}$ the minimal phase $n$th root of
$\phi$,
\begin{align*}
  F_1 &= \Bigl(\tfrac1{\sqrt{n}}\omega^{-k\ell}\Bigr)_{k,\ell=0}^{n-1},
  &
  D_1 &= \diag\bigl(\omega^k\bigr)_{k=0}^{n-1},
\\[1mm]
  F_{\phi} &= \diag\bigl(\phi^{-k/n}\bigr)_{k=0}^{n-1} \cdot  F_1,
  &
  D_{\phi} &= \phi^{1/n} \cdot D_1,
\intertext{and}
  \mathcal{S} &= \Bigl( \sqrt{\tfrac{2}{n+1}} \sin\tfrac{k\ell\pi}{n+1} \Bigr)_{k,\ell=1}^n,
  &
  D_\mathcal{S} &= \diag\bigl( 2\cos\tfrac{k\pi}{n+1} \bigr)_{k=1}^n,
\\[1mm]
  \mathcal{C} &= \Bigl( \sqrt{\tfrac{2}{n}} q_\ell \cos\tfrac{(2k-1)(\ell-1)\pi}{2n} \Bigr)_{k,\ell=1}^n,
  &
  D_\mathcal{C} &= \diag\bigl( 2\cos\tfrac{(k-1)\pi}{n} \bigr)_{k=1}^n.
\end{align*}

The classes listed in Table~\ref{tab:struct} generalize the classical Cauchy,
Toeplitz, Hankel and Vandermonde matrices.
The generators reported below allow one to immerse such structured matrices
into the corresponding displacement class:
\begin{enumerate}
  \item Cauchy matrices $C=(\frac1{t_i-s_j})_{i,j=1}^n$:
    we have $G=H=[1\dots 1]^T$;

  \item Toeplitz matrices $T=(t_{i-j})_{i,j=0}^{n-1}$:
    \begin{equation*}
      G = \left[ \begin{array}{@{}r@{}c@{}lc@{}} 
          & & t_0             & 1 \\
          t_{1-n} & + & t_1             & 0 \\
          & \vdots & & \vdots \\
          t_{-1} & + & t_{n-1}          & 0 
        \end{array} \right],
      \qquad
      \overline{H} = \left[ \begin{array}{@{}cr@{}c@{}l@{}} 
          0      & t_{n-1}&-&t_{-1} \\
          \vdots & & \vdots \\
          0      & t_{1}&-&t_{1-n} \\
          1      & t_{0}
        \end{array} \right];
    \end{equation*}

\item Toeplitz+Hankel matrices $K=(t_{i-j}+h_{i+j})_{i,j=0}^{n-1}$:
  \begin{equation*}
  \small
    \begin{aligned}
      G &= \left[ \begin{array}{@{}r@{}lrrr@{}l@{}} 
          t_0-t_1 & +h_0 & -1 & 0 & t_{-n+1} & +h_{n-1}-h_n \\
          t_1-t_2 & +h_1-h_0 & 0 & 0 & t_{-n+2}-t_{-n+1} & +h_n-h_{n+1} \\
          & \vdots & \vdots & \vdots & & \vdots \\
          t_{n-2}-t_{n-1} & +h_{n-2}-h_{n-3} & 0 & 0 & t_{-1}-t_{-2} & +h_{2n-3}-h_{2n-2} \\
          t_{n-1} & +h_{n-1}-h_{n-2} & 0 & -1 & t_0-t_{-1} & +h_{2n-2}
        \end{array} \right], \\
      H^* &= \left[ \begin{array}{@{}cr@{}lcr@{}lc@{}}
          -1 & & 0 & \cdots & & 0 & 0 \\
          t_{-1} & t_{-2} & +h_0 & \cdots & t_{-n+1} & +h_{n-3} & h_{n-2} \\
          h_n & t_{n-1} & +h_{n+1} & \cdots & t_2 & +h_{2n-2} & t_1 \\
          0 & & 0 & \cdots & & 0 & -1
        \end{array} \right];
    \end{aligned}
  \end{equation*}

\item Vandermonde matrices $W=(w_i^{n-j})_{i,j=1}^n$:
\begin{equation*}
G^T 
= 
\begin{bmatrix} 
  w_1^n-\overline{\phi} & w_2^n-\overline{\phi} & \cdots & w_n^n-\overline{\phi}
\end{bmatrix}, \quad
H^T = \begin{bmatrix} 1 & 0 & \cdots & 0 \end{bmatrix},
\end{equation*}
with $w_i^n\neq\phi$, $i=1,\ldots,n$.
\end{enumerate}
Any Hankel matrix can be trivially transformed into Toeplitz by reverting the
order of its rows, so there is not need to treat this class separately.

Our Matlab routines to convert each displacement structure to Cauchy-like are
listed in Table~\ref{tab:solver}.
For example, \ttt{t2cl} and \ttt{tl2cl} transform a Toeplitz and a
Toeplitz-like matrix, respectively, into a Cauchy-like matrix. A
complete description of each routine in the package is available using
the Matlab \ttt{help} command. Fast routines to apply the matrices
$F_\phi$, $\mathcal{S}$, $\mathcal{C}$, and their inverses, to a vector
are provided in the functions \ttt{ftimes}, \ttt{stimes} and
\ttt{ctimes}; see Table~\ref{tab:aux}.

\begin{remark}
A Cauchy-like matrix is uniquely
identified by its displacement and generators if $t_i\neq s_j$, for any
$i$, $j$. Conversely, if there exist indices $k$, $\ell$ such that
$t_k=s_\ell$, then $\phib_k^* \psib_\ell=0$, and $C_{k\ell}$ cannot be
recovered by~\eqref{eq:Cl}; in such a case the matrix $C$ is said to
be \emph{partially reconstructable}.
\end{remark}

\begin{remark}
A Toeplitz-like matrix of size $m\times n$
has displacement structure with respect to any
pair of displacement matrices $Z_\xi$ and $Z_\eta$, with
$\xi,\eta\in\C$.  The choice $\xi=1$ and
$\eta=\e^{\ii\pi\frac{\gcd(m,n)}m}$ has been proved in~\cite{rod06} to
be optimal, under the constraints $|\xi|=|\eta|=1$, in the sense that it
ensures that the minimum value assumed by the denominator in~\eqref{eq:Cl} is
as large as possible.
\end{remark}

\begin{remark}\label{rem:adjinv}
If the displacement structure of a matrix $A$ is known, it is immediate to
obtain a structured representation for its adjoint and its inverse.
In fact, from~\eqref{eq:EAAF=} it follows that 
$$
\begin{aligned}
F^* A^* - A^* E^* &= -HG^*, \\
F A^{-1} - A^{-1} E &= (-A^{-1} G) (H^* A^{-1}) , \\
\end{aligned}
$$
so that $A^*$ has displacement matrices $(F^*,E^*)$ and generators $(-H,G)$,
while $A^{-1}$ has displacement matrices $(F,E)$.  Its generators
$(\tilde{G},\tilde{H})$ are the solutions of the linear systems
$$
A \tilde{G} = -G, \qquad A^* \tilde{H} = H,
$$
and can be computed by any of the algorithms described in the following
sections.
\end{remark}

Our approach to solve a linear system characterized by any of the above
discussed displacement structures is to preliminarily transform it to
a Cauchy-like system, and to apply an optimized implementation of the
generalized Schur algorithm to the augmented matrix~\eqref{eq:ACb}.
We briefly describe here the displacement structure of this augmented
matrix. Let us consider the linear system~\eqref{eq:Cx=b}, where
$\bm{x},\bm{b}\in\C^{n\times d}$, $d>1$ in the case of multiple right
hand sides. The matrix $\A_{C,\bm{b}}$ associated to the system
inherits a displacement structure, as the following equation holds
\begin{equation}\label{eq:displACB}
  D_{\left[\begin{smallmatrix}\bm{t}\raisebox{1.7mm}{}\\\bm{s}\end{smallmatrix}\right]}
  \,
  \A_{C,\bm{b}} 
  -
  \A_{C,\bm{b}} 
  \, 
  D_{\left[\begin{smallmatrix}\bm{s}\raisebox{1.5mm}{}\\\gamma\bm{e}\end{smallmatrix}\right]}
  =
  \begin{bmatrix} G& (D_{\bm{t}}-\gamma I_n)\bm{b} \\ 0 & 0 \end{bmatrix}
  \begin{bmatrix} H& 0 \\ 0 & I_d \end{bmatrix}^*,
\end{equation}
for any $\gamma\in\C$, where $\bm{e}=[1,\dots,1]^T\in\R^d$ and $G$,
$H$ are given in~\eqref{eq:GHC}. This proves that $\A_{C,\bm{b}}$ is a
Cauchy-like matrix of displacement rank $r+d$.

We assume that $C$ is reconstructable. Even so, the diagonal
entries of the $(2,1)$-block of $\A_{C,\bm{b}}$ are 
nonreconstructable, and the off-diagonal ones are reconstructable if and
only if there are no repetitions in $\bm{s}$, i.e., $s_i\neq s_j$ for
$i\neq j$.
The blocks $(1,2)$ and $(2,2)$ are reconstructable whenever $\gamma\neq t_i,
s_i$, $i=1,\ldots,n$.

It is possible to associate the augmented matrix~\eqref{eq:ACb} to linear
systems having any of the structures reported in Table~\ref{tab:struct}, but
the original structure is preserved only in the Cauchy-like case.
Moreover, these structures are not pivoting invariant.
For these reasons, we do not consider this approach competitive with the one we
follow, in particular for what regards stability.

%
%
\section{Solution of a Cauchy-like linear system}\label{sec:GKO4A}
%
%

Let $C\in\C^{n\times n}$ be a nonsingular Cauchy-like matrix; we want
to solve the linear system~\eqref{eq:Cx=b} where
$\bm{x},\bm{b}\in\C^{n\times d}$. In this Section we describe the core
of our solver, that is an implementation of the generalized Schur
algorithm (GSA) for the augmented Cauchy-like matrix
$\A_{C,\bm{b}}$~\eqref{eq:ACb},
assuming that $C$ has an LU factorization;
several pivoting strategies to treat the general case
will be described in Section~\ref{sec:pivprec}.
We start describing the GSA
for computing the LU factorization of a rectangular Cauchy-like matrix $A$
by recursive Schur complementation, then we apply
such factorization to the case $A=\mathcal{A}_{C,\bm{b}}$.

We introduce some notation which will be
needed throughout this paper. Let $\bm{v}=(v_1,\dots,v_n)^T$ be a
vector of size $n$ and let $A=(a_{ij})$ be an $m\times n$ matrix. 
We use Matlab notation for componentwise division (\ttt{./}) and for
subindexing (\ttt{:}), i.e.,
$$
\begin{aligned}
\bm{v}\ttt{./}\bm{w} &= (v_1/w_1,\ldots,v_n/w_n)^T, \\
A_{2:7,:} &= (a_{ij}), \quad i=2,\dots,7,\ j=1,\dots,n,
\end{aligned}
$$
and define $\hat{\bm{v}}=(v_2,\dots,v_n)^T$ and $\hat{A}=A_{2:m,2:n}$.
Moreover, in the algorithms $a\leftarrow b$ means the usual assignment of $b$
to the variable $a$.

%
%
\subsection{GSA for Cauchy-like matrices}\label{ssec:GKO4C}
%
%
The first step of Gauss algorithm applied to a matrix $A\in\C^{m\times
  n}$ computes the factorization
\begin{equation}\label{eq:gauss1}
  A =
  \begin{bmatrix}
    d & \bm{u}^*\\
    \el & \hat{A}
  \end{bmatrix}
  =
  \begin{bmatrix}
    1 \\
    \frac1d\el & I_{m-1}
  \end{bmatrix}
  \begin{bmatrix}
    d & \bm{u}^*\\
    & \cS_1(A)
  \end{bmatrix},
\end{equation}
where $\cS_1(A)=\hat{A}-\frac1d\el\bm{u}^*$ is the first Schur
complement of $A$. In Lemma~1.1 of~\cite{gko95} it has been proved
that, if $A$ is a Cauchy-like matrix with generators $G\in\C^{m\times
  r}$, $H\in\C^{n\times r}$, and displacement vectors $\bm{t}\in\C^m$,
$\bm{s}\in\C^n$, i.e., if it satisfies the displacement equation
\begin{equation}\label{eq:caudis}
D_\bm{t}A-AD_\bm{s}=GH^*,
\end{equation}
then $\cS_1(A)$ is a Cauchy-like matrix too, satisfying the equation
\begin{equation*}
  D_{\hat{\bm{t}}}\,\cS_1(A)-\cS_1(A)\,D_{\hat{\bm{s}}} = 
  \tilde{G}\,\tilde{H}^*.
\end{equation*}
Here the generators $\tilde{G}$, $\tilde{H}$ are defined by
\begin{equation}\label{eq:GH+}
  \begin{bmatrix}0\\\tilde{G}\end{bmatrix}
  =
  G
  -
  \begin{bmatrix}1\\\frac1d\el\end{bmatrix}G_{1,:}\,,
  \qquad
  \begin{bmatrix}0\\\tilde{H}\end{bmatrix}
  =
  H
  -
  \begin{bmatrix}1\\\frac{1}{\overline{d}}\bm{u}\end{bmatrix}H_{1,:}\,.
\end{equation}

The GSA consists of the recursive application 
of the decomposition~\eqref{eq:gauss1} to $\cS_1(A)$,
operating on the data which define the displacement
structure. At the $k$-th iteration, the quantities $d$, $\el$ and
$\bm{u}$ are computed using formula~\eqref{eq:Cl};
then the generators are updated according to~\eqref{eq:GH+}. It's important to
remark that the step $k$ of GSA computes the generators of $\cS_k(A)$, since
$\cS_k(A)=\cS_1\bigl(\cS_{k-1}(A)\bigr)$, and that each $\cS_k(A)$ is
reconstructable if $A$ is.

If the LU factors of $A$ are required, then the vectors $\frac1d\el$
and $[d\;\bm{u}^*]$ must be stored at each iteration: after the $p$-th
iteration, $p(m+n-p)$ memory locations are required to store these vectors, and a
partial LU factorization of $A$ is computed
\begin{equation}\label{eq:fCl}
  A
  =
  \begin{bmatrix}
    A_{11} & A_{12} \\
    A_{21} & A_{22}
  \end{bmatrix}
  =
  \begin{bmatrix}
    L &    \\
    A_{21}U^{-1} & I_{m-p}
  \end{bmatrix}
  \begin{bmatrix}
    U & L^{-1}A_{12}    \\
      & \cS_p(A)
  \end{bmatrix},
\end{equation}
where $A_{11}=LU\in\C^{p\times p}$. If the LU factors are
not required, it is possible to discard $d$, $\el$ and $\bm{u}$ after
the computation of $\tilde{G}$ and $\tilde{H}$. This is the case we
are interested in, since our aim is to compute $\cS_n(\A_{C,\bm{b}})$.

The GSA for the computation of the Schur complement $\cS_p(A)$ is described in
Algorithm~\ref{alg:GKOminimal}. The memory locations required by the
algorithm are $(m+n)(r+2)$. The computational cost is
$p((4r+1)\alpha+1)$ floating point operations (\flops), if the input is
real, and $2p(5(\alpha+1)+8r(\alpha+\frac14))$ \flops, if the input is
complex, where $\alpha=m+n-p$, and assuming 2 \flops\ for each complex
sum, 6 for each product and 10 for each division. The algorithm
outputs the displacement matrices and the generators of $\cS_p(A)$.
\begin{algorithm}[ht]
  \begin{algo}
    \REQUIRE $m,n,r,p \in\N^+$ s.t. $p<\min(m,n)$
    \REQUIRE $\bm{t}\in\C^m$, $\bm{s}\in\C^n$, $G\in\C^{m\times r}$,
    $H\in\C^{n\times r}$ s.t. $t_i\neq s_j$, $\forall i,j$
    \STATE   $H\leftarrow H^*$
    \FOR     $k=1$ \To $p$
    \STATE   $d\leftarrow (G_{k,:}\cdot H_{:,k})/(t_k-s_k)$
    \STATE   \If $d=0$ \Then \textsc{stop}
    \STATE   $\el\leftarrow (G_{k+1:m,:}\cdot H_{:,k})./(t_{k+1:m}-s_k)$
    \STATE   $\bm{u}\leftarrow (G_{k,:}\cdot H_{:,k+1:n})./(t_k-s_{k+1:n}^T)$
    \STATE   $G_{k+1:m,:}\leftarrow G_{k+1:m,:}-\el\cdot(\frac1d G_{k,:})$
    \STATE   $H_{:,k+1:n}\leftarrow H_{:,k+1:n}-(\frac1d H_{:,k})\cdot\bm{u}$
    \ENDFOR
    \OUTPUT  $\bm{t}_{p+1:m}$, $\bm{s}_{p+1:n}$, $G_{p+1:m,:}$, $(H_{:,p+1:n})^*$
  \end{algo}
  \caption{GSA for computing $\cS_p(A)$, assuming
  $D_{\bm{t}}A-AD_{\bm{s}}=GH^*$.}
  \label{alg:GKOminimal}
\end{algorithm}

\begin{remark}\label{rem:ts}
  Algorithm~\ref{alg:GKOminimal} requires that the Cauchy-like matrix
  $A$ is reconstruct\-ab\-le. If it is not, the algorithm must be
  modified in order to store the nonreconstructable entries. 
  This can be done by employing part of the
  generators, since after the step $k$ the first $k$ rows of both $G$
  and $H$ become unused.
\end{remark}
%
%
\subsection{GSA for the augmented matrix $\A_{C,\bm{b}}$}\label{ssec:GKO4A}
%
%
To compute $\cS_n(\A_{C,\bm{b}})$ we apply $n$ steps of the GSA to the
matrix $\A_{C,\bm{b}}$; in this case the factorization~\eqref{eq:fCl}
reads
\begin{equation}\label{eq:fmb}
	\A_{C,\bm{b}}
	=
	\begin{bmatrix}
		C  & \bm{b} \\
		-I_n & 0
	\end{bmatrix}
	=
	\begin{bmatrix}
		L &    \\
		-U^{-1} & I_n
	\end{bmatrix}
	\begin{bmatrix}
		U & L^{-1}\bm{b}    \\
	      & \cS_n(\A_{C,\bm{b}})
	\end{bmatrix},
\end{equation}
and it emphasizes the fact that $L$ and $U^{-1}$ are computed by columns, and
$U$ by rows. We observe that, in this case,
the reconstructability requirement in line 2 of Algorithm~\ref{alg:GKOminimal} 
is not satisfied. This
causes no problems when $s_i\neq s_j$, $\forall i\neq j$, since the
nonreconstructable entries of $\A_{C,\bm{b}}$
are located along the diagonal of its
$(2,1)$-block, and are known to be $-1$. On the contrary, in the
presence of repeated entries in the vector $\bm{s}$,
Algorithm~\ref{alg:GKOminimal} is not applicable; an
algorithm to treat this case will be proposed in Section~\ref{ssec:ss}.

One can simultaneously compute the solution of $d$ linear systems, by setting
$\bm{x},\bm{b}\in\C^{n\times d}$ in Algorithm~\ref{alg:GKOminimal};
see~\eqref{eq:displACB}. 
A straightforward application yields a cost of $(8r+8d+2)n^2+O(n)$ \flops\ if 
the input is real, and $(32r+32d+20)n^2+O(n)$ \flops\ if the input is
complex. This cost can be reduced thanks to the following remarks.

Since in general the right hand side $\bm{b}$ is unstructured, it is 
convenient to store it as a
matrix, rather than by means of the displacement relation~\eqref{eq:displACB}.
We adapted our algorithm in order to explicitly
store the $d$ rightmost columns of $\A_{C,\bm{b}}$,
i.e., $\left[\begin{smallmatrix}\bm{b}\\0\end{smallmatrix}\right]$, and
to perform ``traditional'' Gauss elimination on them. Accordingly, we
used the displacement equation
\begin{equation}\label{eq:displAS}
  \cD_{\left[\begin{smallmatrix}\bm{t}\raisebox{1.7mm}{}\\\bm{s}\end{smallmatrix}\right]}
  \, 
  \left[\begin{array}{@{}c@{}}C\\-I_n\end{array}\right]
  - \left[\begin{array}{@{}c@{}}C\\-I_n\end{array}\right] \, D_{\bm{s}}
  =
  \begin{bmatrix} G_C \\ 0  \end{bmatrix}
  H_C^*
\end{equation}
to represent the $n$ leftmost columns of $\A_{C,\bm{b}}$.

By the particular structure of
$\A_{C,\bm{b}}$ (see~\eqref{eq:fmb}), the step $k$ of
Algorithm~\ref{alg:GKOminimal} yields a
vector $\el\in\C^{2n-k}$ whose last $n-k$ entries vanish, so that it
can be stored in a vector of length $n$. A consequence is that, at step $k$,
the GSA only needs to modify the rows from $k+1$ to $n+k$ of
$\A_{C,\bm{b}}$, and to do so only such rows are changed in the left generator.
Given its initial form (see~\eqref{eq:displAS}), only
$n$ rows of the left generator are essential, regardless the index $k$. To
optimize the memory access, we reuse the first $k$ rows of $G_C$ to
store the rows of the left generator ranging from $n+1$ to
$n+k$; see also Remark~\ref{rem:ts}. The same procedure is applied
to the vector $\el$ and to the submatrix
$\left[\begin{smallmatrix}\bm{b}\\0\end{smallmatrix}\right]$, whose
essential part is stored in $\bm{b}$.

\begin{algorithm}[ht]
  \begin{algo}
    \REQUIRE $n,r,d \in\N^+$
    \REQUIRE $\bm{t},\bm{s}\in\C^n$; $G,H\in\C^{n\times r}$; $\bm{b}\in\C^{n\times d}$
    \STATE[] \phantom{\algofont{Require: }}s.t. $t_i\neq s_j$, $\forall i,j$, 
    and $s_i\neq s_j$, $\forall i\neq j$
    \STATE $H\leftarrow H^*$
    \FOR   $k=1$ \To $n$
    \STATE $\el\leftarrow(G\cdot H_{:,k})./([s_{1:k-1};t_{k:n}]-s_k)$
    \STATE \If $\ell_k=0$ \Then \textsc{stop}
    \STATE $\bm{u}\leftarrow(G_{k,:}\cdot H_{:,k+1:n})./(t_k-s_{k+1:n}^T)$
    \STATE $d\leftarrow \ell_k$
    \STATE$\ell_k\leftarrow -1$
    \STATE $\bm{g}\leftarrow\frac1dG_{k,:}$
    \STATE $G_{k,:}\leftarrow\bm{0}$
    \STATE $G\leftarrow G-\el\cdot\bm{g}$
    \STATE $\bm{f}\leftarrow\frac 1 d \bm{b}_{k,:}$
    \STATE $\bm{b}_{k,:}\leftarrow\bm{0}$
    \STATE $\bm{b}\leftarrow\bm{b} -\el\cdot\bm{f}$
    \STATE $H_{:,k+1:n}\leftarrow H_{:,k+1:n}-(\frac1d H_{:,k})\cdot\bm{u}$
    \ENDFOR
    \OUTPUT $\bm{x}\leftarrow\bm{b}$
  \end{algo}
  \caption{Solver for $C\bm{x}=\bm{b}$, assuming
  $D_{\bm{t}}C-CD_{\bm{s}}=GH^*$ ($s_i\neq s_j$).}
  \label{alg:GKOstructured}
\end{algorithm}
This leads to Algorithm~\ref{alg:GKOstructured}.
We observe that it approximately needs $(2r+d+2)n$ memory locations: the
storage is then linear in $n$, while the original GKO algorithm~\cite{gko95}
requires $O(n^2)$ locations. Moreover, the storage required by
Algorithm~\ref{alg:GKOstructured} is almost minimal, since it uses only two
extra $n$-vectors besides the right hand side $\bm{b}$ and the
displacement data of $C$. The required \flops\ are
$(6r+2d+\frac32)n^2+O(n)$ in the real case and $(24r+8d+15)n^2+O(n)$
in the complex case. If we compare it to Algorithm~\ref{alg:GKOminimal} when
$r=2$, which is the case that occurs when we
solve Toeplitz systems, the required \flops\ scale to $61\%$ if $d=1$
and $45\%$ if $d=4$.

The requirements $t_i\neq s_j$, $\forall i,j$, and $s_i\neq s_j$,
$\forall i\neq j$, guarantee that $\A_{C,\bm{b}}$ is
reconstructable, apart from the diagonal entries of its
$(2,1)$-block.

%
%
\subsection{The case of multiple $s_i$}\label{ssec:ss}
%
%
As already observed in Section~\ref{ssec:GKO4A}, if some
repetitions occur in the vector $\bm{s}$,
then Algorithm~\ref{alg:GKOstructured}
is not applicable, because the upper triangular part of the
$(2,1)$-block of $\A_{C,\bm{b}}$ has some off-diagonal entries which
are not reconstructable, and thus the vector $\el$ cannot be computed
by line \ttt{5} of Algorithm~\ref{alg:GKOstructured}. 
It is worth noting that no entry of $\bm{s}$
(or of $\bm{t}$) can be repeated more than $r$ times, otherwise $C$ would be
singular.

To overcome this difficulty there are at least two possibilities. Both are
based on a permutation of $\bm{s}$, which induces a
change of variable in $C\bm{x}=\bm{b}$, as shown in the following trivial
proposition.

%
%

\begin{proposition}\label{prop:ss}
  Let $C\in\C^{n\times n}$, $\bm{t},\bm{s}\in\C^n$, $G,H\in\C^{n\times
    r}$, $\bm{b}\in\C^{n\times d}$, and assume that
  \begin{equation*}
    D_{\bm{t}}C-CD_{\bm{s}}=GH^*.
  \end{equation*}
  If $Q$ is a permutation matrix, then
  \begin{equation*}
    D_{\bm{t}}\widetilde{C}-\widetilde{C}D_{\tilde{\bm{s}}}=G\widetilde{H}^*,
  \end{equation*}
  where $\widetilde{C}=CQ^T$, $\tilde{\bm{s}}=Q\bm{s}$, and $\widetilde{H}=QH$. 
  When $C$ is nonsingular, the solution of
  $C\bm{x}=\bm{b}$ is related to the one of
  $\widetilde{C}\tilde{\bm{x}}=\bm{b}$ via $\tilde{\bm{x}}=Q\bm{x}$.
\end{proposition}

%
%
\subsubsection{Redundant-injective splitting of\/ $\bm{s}$}\label{ssec:ss1}
%
%
The simplest way to deal with repeated entries in $\bm{s}$ in
Algorithm~\ref{alg:GKOstructured} is to choose $Q$ such that
$\tilde{\bm{s}}=Q\bm{s}=\left[\begin{smallmatrix}\tilde{\bm{s}}_1\\\tilde{\bm{s}}_2\end{smallmatrix}\right]$,
where $\tilde{\bm{s}}_2\in\C^{\nu}$ has no repetitions and $\nu$ is as
large as possible. 
We also partition
$\tilde{\bm{x}}=Q\bm{x}=\left[\begin{smallmatrix}\tilde{\bm{x}}_1\\\tilde{\bm{x}}_2\end{smallmatrix}\right]$
and
$\bm{b}=\left[\begin{smallmatrix}\bm{b}_1\\\bm{b}_2\end{smallmatrix}\right]$
so that $\tilde{\bm{x}}_2,\bm{b}_2\in\C^{\nu\times d}$.
Since $C$ is invertible, $\nu\geqslant\lceil n/r\rceil$.

By adapting Algorithm~\ref{alg:GKOstructured} to the augmented matrix
\begin{equation*}
  \B_{\widetilde{C},\bm{b}}
  =
  \begin{bmatrix}
    \widetilde{C}&\bm{b}\\ 
    J_\nu&0
  \end{bmatrix},
  \qquad J_\nu=\begin{bmatrix}0&-I_\nu\end{bmatrix}
  \in\R^{\nu\times n},
\end{equation*}
we obtain $\cS_n(\B_{\widetilde{C},\bm{b}})=-J_\nu
\widetilde{C}^{-1}\bm{b}=\tilde{\bm{x}}_2$. 
Then, the vector $\tilde{\bm{x}}_1$ can be computed by solving the system
\begin{equation}\label{brutto}
  \widetilde{C}_{11}\tilde{\bm{x}}_{11} =
  (\bm{b}_{1}-\widetilde{C}_{12}\tilde{\bm{x}}_2),
\end{equation}
where the matrix
$\bigl[\widetilde{C}_{11}\,\widetilde{C}_{12}\bigr]$
contains the first $n-\nu$ rows of $\widetilde{C}$.

 

This approach is simple, but has two main disadvantages. The first one
is that we need extra storage for $\widetilde{C}_{11}$ and
$\widetilde{C}_{12}$, since the initial generators get overwritten
during the algorithm. The second is that it is reasonable to solve
the system~\eqref{brutto} by an unstructured method only if $n-\nu$ is
small. If this condition is not met, a possible strategy is to apply
recursively the same technique to $\widetilde{C}_{11}$, but this approach 
would lead, in the worst case, to $r-1$ recursive steps.

\begin{remark}
  We observe that when applying our Matlab implementation of
  Algorithm~\ref{alg:GKOstructured} to a Cauchy-like matrix with
  repetitions in $\bm{s}$, only $n-\nu$ components of the computed solution
  vector are affected by the overflows caused by nonreconstructable
  entries (i.e., they are\/ \verb|Inf| or\/ \verb|NaN|), while the remaining
  $\nu$ components are correct.
  This is due to Matlab full implementation of IEEE floating-point arithmetic.
\end{remark}

%
%
\subsubsection{Gathering of\/ $\bm{s}$}\label{ssec:ss2}
%
%
A different way to deal with repeated entries in $\bm{s}$
happens to be very effective both from the point of view of storage
and computational cost.

Let $\sigma_j$, $j=1,\dots,\nu$, be the distinct entries of $\bm{s}$,
each with multiplicity $\mu_j$, and choose a permutation
matrix $Q$ such that each subset of repeated components is gathered, i.e.,
occupies contiguous entries in
$\tilde{\bm{s}}=Q\bm{s}$. It follows that $\tilde{\bm{s}}$ can be
partitioned as
\begin{equation*}
  \tilde{\bm{s}}=\begin{bmatrix}\sigmab^{(1)}\\[-1mm]\vdots\\[1mm]\sigmab^{(\nu)}\end{bmatrix},
  \quad\textnormal{where}\quad
  \sigmab^{(j)}=\sigma_j\left[\begin{smallmatrix}1\\[-1mm]\vdots\\[1mm]1\end{smallmatrix}\right]\in\C^{\mu_j}
  \quad\textnormal{and}\quad
  \sum_{j=1}^{\nu}\mu_j=n.
\end{equation*}
The order of the $\sigma_j$ is not important and $Q$ is not unique.
Since $C$ is nonsingular, $\mu_j\leqslant r$, $j=1,\dots,\nu$.

Changing variable according to Proposition~\ref{prop:ss} causes the
nonreconstructable entries in the $(2,1)$-block of
$\A_{\widetilde{C},\bm{b}}$ to be grouped in the square blocks
$\mathcal{T}_j$, $j=1,\dots,\nu$, each one of size $\mu_j$, located
along the diagonal of the $(2,1)$-block. If a $\sigma_j$ occurs only
once in $\tilde{\bm{s}}$ ($\mu_j=1$) then $\mathcal{T}_j$
is a scalar.

If we apply Algorithm~\ref{alg:GKOstructured} to the system
$\widetilde{C}\tilde{\bm{x}}=\bm{b}$, the vector $\el$ computed at
step $k$ intersects only one nonreconstructable block, say
$\mathcal{T}_{j(k)}$, and we define $\alpha_k$ to be the column index
in $\A_{\widetilde{C},\bm{b}}$ of the first column of
$\mathcal{T}_{j(k)}$. Since the $(2,1)$-block of
$\A_{\widetilde{C},\bm{b}}$ is upper triangular and has
$-1$ along its diagonal, during the algorithm we need to explicitly
store and update the strictly upper triangular part of each
$\mathcal{T}_j$. Moreover, at step $k$ we need to store only the first
$k-\alpha_k+1$ rows of $\mathcal{T}_{j(k)}$, i.e., those ranging from
$n+\alpha_k$ to $n+k$.  These rows fit in a natural way into
$H_{\alpha_k:k,:}$ since $\mu_{j(k)}\leqslant r$; see also Remark~\ref{rem:ts}.
The update of $\mathcal{T}_{j(k)}$ is performed by
standard Gaussian elimination, as we do for the rightmost columns of
$\A_{\widetilde{C},\bm{b}}$.

The above technique can be embedded in Algorithm~\ref{alg:GKOstructured}, to
make it applicable in the presence of
repeated entries in the vector $\bm{s}$. The result is a general
method, outlined in Algorithm~\ref{alg:GKOss}, for the solution of a
Cauchy-like linear system whose coefficient matrix is reconstructable
and nonsingular. Besides the vector
$\alphab=[\alpha_1,\dots,\alpha_n]^T$, we use the auxiliary vector
$\omegab=[\omega_1,\dots,\omega_n]^T$, where $\omega_k$ is the column
index in $\A_{\widetilde{C},\bm{b}}$ of the last column of
$\mathcal{T}_{j(k)}$. The vectors $\alphab$ and $\omegab$ can be
constructed in Matlab with library functions at negligible cost;
obviously, in the code the computation involving the permutation
matrix $Q$ is handled via a vector of indices.

\begin{algorithm}
  \begin{algo}
    \REQUIRE      $n,r,d \in\N^+$
    \REQUIRE      $\bm{t},\bm{s}\in\C^n$; $G,H\in\C^{n\times r}$; $\bm{b}\in\C^{n\times d}$ s.t. 
                  $t_i\neq s_j$, $\forall i,j$
    \COMPUTE[3.1:]$Q$, $\alphab$, $\omegab$ (see text)
    \STATE[3.2:]  $H\leftarrow (QH)^*$
    \STATE[3.3:]  $\bm{s}\leftarrow Q\bm{s}$
\algosetline{4}
    \FOR          $k=1:n$
    \STATE[5.1:]  $\el_{1:\alpha_k-1}\leftarrow(G_{1:\alpha_k-1,:}\cdot H_{:,k})./(s_{1:\alpha_k-1}-s_k)$
    \STATE[5.2:] \If$\alpha_k<k$
    \Then  $\el_{\alpha_k:k-1}\leftarrow (H_{k-\alpha_k,\alpha_k:k-1})^T$
    \STATE[5.3:]  $\el_{k:n}\leftarrow(G_{k:n,:}\cdot H_{:,k})./(t_{k:n}-s_k)$
    \STATE[\dots] \emph{\{ lines from 6 to 16 of Algorithm~\ref{alg:GKOstructured} \}}
    \IF[16.1:] $k<\omega_k$
    \STATE[16.2:] $\gamma\leftarrow k-\alpha_k+1$
    \STATE[16.3:] $\delta\leftarrow \omega_k-\alpha_k$
    \STATE[16.4:] $H_{\gamma:\delta,k}\leftarrow 0$
    \STATE[16.5:] $H_{\gamma:\delta,\alpha_k:k}\leftarrow H_{\gamma:\delta,\alpha_k:k} - (\frac1d\bm{u}_{1:\omega_k-k})^T\cdot(\el_{\alpha_k:k})^T$
    \ENDIF[16.6:]
\algosetline{17}
    \ENDFOR
    \OUTPUT       $\bm{x}\leftarrow Q^T\bm{b}$
  \end{algo}
  \caption{Solver for $C\bm{x}=\bm{b}$, assuming
  $D_{\bm{t}}C-CD_{\bm{s}}=GH^*$ (mult.~knots).}
  \label{alg:GKOss}
\end{algorithm}

Whenever the vector $\el$ intersects a block $\mathcal{T}_j$, some of its
entries have to be extracted from $H$; see lines
\ttt{5.1}--\ttt{5.3} in Algorithm~\ref{alg:GKOss}.
The update of the nonreconstructable
block $\mathcal{T}_{j(k)}$ is performed at lines
\ttt{16.1}--\ttt{16.6}. This strategy of handling
repetitions in $\bm{s}$ does not increase the
overall complexity with respect to Algorithm~\ref{alg:GKOstructured}.

We note that this approach can be numerically effective 
also when some entries of the displacement vector $\bm{s}$ are very close.
In this case, in order to avoid small denominators in formula~\eqref{eq:Cl}, it
may be convenient to approximate the system matrix
by collapsing the clustered entries in $\bm{s}$, and then apply
Algorithm~\ref{alg:GKOss}.
This possibility will be explored in Section~\ref{sec:numres}.

%
%
\section{Pivoting strategies}
\label{sec:pivprec}
%
%
Algorithms~\ref{alg:GKOstructured} and~\ref{alg:GKOss} require
that the LU factorization of $C$ exists. Both to avoid this requirement,
and for stability issues, it is important to introduce a pivoting strategy
in the algorithms. Any such strategy should not increase the computational
cost and the amount of memory required, which are $O(n^2)$ and $O(n)$,
respectively. This condition is met by partial pivoting, while 
complete pivoting would raise the computational cost to $O(n^3)$, since it
requires the explicit computation of a full matrix at each step. In~\cite{sb95}
and~\cite{gu98}, the authors propose some approximations to
complete pivoting, which keep the computational cost 
quadratic. We adapted these strategies to our
augmented matrix approach for solving a Cauchy-like linear system, 
keeping the storage linear. We
implemented three variations of Algorithm~\ref{alg:GKOstructured},
in the case there are no repetitions in $\bm{s}$, to include:
\begin{itemize}
  \item[a)] partial pivoting (Algorithm~\ref{alg:GKOssPP});
  \item[b)] Sweet \& Brent's pivoting~\cite{sb95} (Algorithm~\ref{alg:GKO_SB});
  \item[c)] Gu's pivoting and generator scaling technique~\cite{gu98} (Algorithm~\ref{alg:GKO_GU}).
\end{itemize}

We also implemented partial pivoting for Algorithm~\ref{alg:GKOss}, i.e., for
the case when there are repetitions in $\bm{s}$; see
Algorithm~\ref{alg:GKOssPP}. 
In this case, only partial pivoting is suitable, since it preserves the
ordering of $\tilde{\bm{s}}$,
as it will be made clear in the following.
For comparison issues only, we included complete pivoting in
Algorithm~\ref{alg:GKO_TP}.



Since our goal is to compute the solution to system~\eqref{eq:Cx=b} as the
Schur complement of order $n$ of the augmented matrix~\eqref{eq:ACb}, pivoting
must be performed only on the first $n$ rows and columns of $\A_{C,\bm{b}}$
(or $\A_{\widetilde{C},\bm{b}}$).
If we operate only on rows, factorization~\eqref{eq:fmb} is replaced by
\begin{equation}\label{eq:fPmb}
  \mathcal{P}
  \A_{C,\bm{b}}
  =
  \begin{bmatrix}
    P &  \\
    & I_n  
  \end{bmatrix}
  \begin{bmatrix}
    C & \bm{b} \\
    -I_n  & 0
  \end{bmatrix}
  =
  \begin{bmatrix}
    L       &   \\
    -U^{-1} & I_n
  \end{bmatrix}
  \begin{bmatrix}
    U & L^{-1}P\bm{b} \\
      & C^{-1}\bm{b} 
  \end{bmatrix},
\end{equation}
while if both rows and columns are permuted, it is convenient to apply the
factorization
\begin{equation}\label{eq:genpiv}
\begin{aligned}
  \mathcal{P}
  \A_{C,\bm{b}}
  \mathcal{Q}^T
  &=
  \begin{bmatrix}
    P &  \\
      & Q
  \end{bmatrix} 
  \begin{bmatrix}
    C    & \bm{b} \\
    -I_n & 0
  \end{bmatrix}
  \begin{bmatrix}
    Q^T &  \\
      & I_d 
  \end{bmatrix} \\
  &=
  \begin{bmatrix}
    L       &   \\
    -U^{-1} & I_n
  \end{bmatrix}
  \begin{bmatrix}
    U & L^{-1}P\bm{b} \\
      & QC^{-1}\bm{b}                  
  \end{bmatrix}.
\end{aligned}
\end{equation}
In~\eqref{eq:genpiv}, when columns $i$ and $j$ are swapped, we also
exchange rows $n+i$ and $n+j$ of the augmented matrix.
This unusual permutation does not affect the computation, as the GSA performs
the same arithmetic operations using a different rows ordering, with the only
difference of computing a permutation $Q\bm{x}$ of the system solution.
The advantage of this approach is that the nonreconstructable entries stay
along the diagonal of the $(2,1)$-block.

As mentioned above, we applied only partial pivoting to
Algorithm~\ref{alg:GKOss}, as column pivoting induces a permutation in the
right displacement vector $\tilde{\bm{s}}$.
So, in this case only factorization~\eqref{eq:fPmb} is employed, replacing 
$C$ by $\widetilde{C}$.

%
%
\subsection{Partial pivoting}\label{ssec:pp}
%
%
Partial pivoting can be included in Algorithms~\ref{alg:GKOstructured}
and~\ref{alg:GKOss} by inserting two instructions before line
\ttt{6}, as shown in Algorithm~\ref{alg:GKOssPP}.
\begin{algorithm}[ht]
  \begin{algo}
\algoskipindent{1}
    \STATE[\dots] \emph{\{ Insert in Algorithm~\ref{alg:GKOstructured} and~\ref{alg:GKOss}: \}}
    \STATE[5.4:]  \textsc{find} $i$ \textsc{such that} $|\ell_i|=\max|\el_{k:n}|$
    \STATE[5.5:]     \If $i\neq k$ \Then
        \textsc{swap rows} $k,i$ \textsc{of} $G$, $\bm{t}$, $\el$, $\bm{b}$
    \STATE[\dots]
  \end{algo}
  \caption{Solver for $C\bm{x}=\bm{b}$ with partial pivoting.}
  \label{alg:GKOssPP}
\end{algorithm}

%
%
%
\subsection{Sweet \& Brent pivoting}
%
%

In~\cite{sb95}, an error analysis is performed for the LU factorization computed
by the GSA with partial pivoting, in the case of Cauchy and Toeplitz matrices.
Sweet and Brent obtained an upper bound for the error depending on the LU
factors, as one might expect, and on the \emph{generator growth factors}.
In fact they showed that, when the displacement rank is larger than 1, there
can be a large growth in the generators entries, and this is reflected in the
solution error.

In order to approximate complete pivoting, which would increase the stability
of the algorithm, Sweet and Brent proposed to choose the pivot, at step $k$, by
searching both the $k$th column and the $k$th row.
The maxima
\begin{equation*}
  p_1=|\ell_{i_1}|=\max_{i=k+1,\ldots,n}|\ell_i|
  \qquad
  \textnormal{and}
  \qquad
  p_2=|u_{i_2}|=\max_{j=k+1,\ldots,n}|u_j|,
\end{equation*}
are compared with the actual pivot $d$; see~\eqref{eq:gauss1}.
If $d<\max\{p_1,p_2\}$, then if $p_1\geqslant p_2$ the rows $k$ and $i_1$ are
swapped, otherwise the columns $k$ and $i_2$ are swapped.
Applying this pivoting strategy to the augmented matrix~\eqref{eq:ACb},
whenever $d<p_1<p_2$ we also swap the rows $n+k$ and $n+i_2$, according to
factorization~\eqref{eq:genpiv}. This originates Algorithm~\ref{alg:GKO_SB}.

\begin{algorithm}[ht]
  \begin{algo}
    \REQUIRE $n,r,d \in\N^+$
    \REQUIRE $\bm{t},\bm{s}\in\C^n$; $G,H\in\C^{n\times r}$; $\bm{b}\in\C^{n\times d}$
    \STATE[] \phantom{\algofont{Require: }}s.t. $t_i\neq s_j$, $\forall i,j$, 
    and $s_i\neq s_j$, $\forall i\neq j$
    \STATE $H\leftarrow H^*$
    \STATE[3.9:]  $Q \leftarrow I_n$
    \FOR   $k=1$ \To $n$
    \STATE $\el_{k:n}\leftarrow(G_{k:n}\cdot H_{:,k})./(t_{k:n}-s_k)$
    \STATE[5.1:]  $\bm{u}_{k+1:n}\leftarrow(G_{k,:}\cdot H_{:,k+1:n})./(t_k-s_{k+1:n}^T)$
    \STATE[5.2:]  \textsc{find} $i_1$ \textsc{such that} $p_1:=|\ell_{i_1}|=\max|\el_{k:n}|$
    \STATE[5.3:]  \textsc{find} $i_2$ \textsc{such that} $p_2:=|u_{i_2}|=\max|\bm{u}_{k+1:n}|$
    \IF       $p_1=p_2=0$ \Then \textsc{stop}
    \ELSEIF[6.1:] $p_2>p_1$
    \STATE[6.2:] \textsc{swap columns} $k,i_2$ \textsc{of} $H$, $\bm{s}$, $Q^T$
    \STATE[6.3:] $u_{i_2} \leftarrow \ell_k$
    \STATE[6.4:] $\el_{k:n}\leftarrow(G_{k:n}\cdot H_{:,k})./(t_{k:n}-s_k)$
    \ELSEIF[6.5:] $i_1>k$
    \STATE[6.6:] \textsc{swap rows} $k,i_1$ \textsc{of} $G$, $\bm{t}$, $\el$, $\bm{b}$
    \STATE[6.7:]  $\bm{u}_{k+1:n}\leftarrow(G_{k,:}\cdot H_{:,k+1:n})./(t_k-s_{k+1:n}^T)$
    \ENDIF[6.8:]
    \STATE[7:] $\el_{1:k-1}\leftarrow(G_{1:k-1}\cdot H_{:,k})./(s_{1:k-1}-s_k)$
    \STATE[\dots]  \emph{\{ Insert lines 8--16 from Algorithm~\ref{alg:GKOstructured} \}}
    \algosetline{17}
    \ENDFOR
    \OUTPUT       $\bm{x}\leftarrow Q^T\bm{b}$
  \end{algo}
  \caption{Solver for $C\bm{x}=\bm{b}$ with Sweet and Brent's pivoting.}
  \label{alg:GKO_SB}
\end{algorithm}

%
%
\subsection{Gu's pivoting}
%
%

In~\cite{gu98} a different approach was proposed to prevent the generator
growth.
At each step a compact QR factorization of the left generator $G$
(see~\eqref{eq:caudis}) is computed, and the triangular factor is
\emph{transferred} to the right generator, i.e.,
\begin{equation*}
  GH^* = (UR)H^* = U (HR^*)^* = \tilde{G}\tilde{H}^*,
\end{equation*}
where $R\in\C^{r\times r}$ is upper triangular, and $U$ is a matrix with
orthonormal columns having the same size than $G$.
This procedure has the effect of keeping the 2-norm of the left generator
constant across the iterations.
Moreover, it is immediate to observe that the $j$th column of the right
generator $\tilde{H}^*$ has the same 2-norm than the corresponding column of
the product $GH^*$.

Gu's approximation of complete pivoting consists of selecting the column
$j_\text{max}$ of $\tilde{H}^*$ having the largest 2-norm, and swapping
the columns $k$ and $j_\text{max}$ of the system matrix before performing
partial pivoting.
Moreover, he observes that in most cases it is sufficient to apply this
procedure every $K$ steps, instead than at each iteration $k$, and uses the
value $K=10$ in his numerical experiments.

To adapt Gu's technique to our augmented matrix approach, we compute the QR
factorization of the last $n-k+1$ rows of $G$, and update the generators
accordingly; see Algorithm~\ref{alg:GKO_GU}.

\begin{algorithm}[ht]
  \begin{algo}
    \REQUIRE $n,r,d,K \in\N^+$
    \REQUIRE $\bm{t},\bm{s}\in\C^n$; $G,H\in\C^{n\times r}$; $\bm{b}\in\C^{n\times d}$
    \STATE[] \phantom{\algofont{Require: }}s.t. $t_i\neq s_j$, $\forall i,j$, 
    and $s_i\neq s_j$, $\forall i\neq j$
    \STATE $H\leftarrow H^*$
    \STATE[3.9:]  $Q \leftarrow I_n$
    \FOR   $k=1$ \To $n$
    \IF[4.1:] mod$(k,K)=1$ \textsc{and} $k\leqslant n-r+1$
    \STATE[4.2:] $[U,\,R] = \textnormal{qr}(G_{k:n,:},0)$ \COMMENT compact QR factorization
    \STATE[4.3:] $G_{1:k-1,:} \leftarrow G_{1:k-1,:}R^{-1}$
    \STATE[4.4:] $G_{k:n,:} \leftarrow U$
    \STATE[4.5:] $H_{:,k:n}\leftarrow R\cdot H_{:,k:n}$
    \STATE[4.6:] \textsc{find} $i_2$ \textsc{such that} $||H_{:,i_2}||_2=\max_{k\leqslant j \leqslant n}||H_{:,j}||_2$
    \STATE[4.7:] \If $i_2\neq k$ \Then \textsc{swap columns} $k,i_2$ \textsc{of} $H$, $\bm{s}$, $Q^T$
    \ENDIF[4.8:]
    \STATE $\el\leftarrow(G\cdot H_{:,k})./([s_{1:k-1};t_{k:n}]-s_k)$
    \STATE[5.4:]  \textsc{find} $i_1$ \textsc{such that} $|\ell_i|=\max|\el_{k:n}|$
    \STATE[5.5:]     \If $i_1\neq k$ \Then
        \textsc{swap rows} $k,i_1$ \textsc{of} $G$, $\bm{t}$, $\el$, $\bm{b}$
    \STATE[\dots]  \emph{\{ Insert lines 6--16 from Algorithm~\ref{alg:GKOstructured} \}}
    \algosetline{17}
    \ENDFOR
    \OUTPUT       $\bm{x}\leftarrow Q^T\bm{b}$
  \end{algo}
  \caption{Solver for $C\bm{x}=\bm{b}$ with Gu's pivoting.}
  \label{alg:GKO_GU}
\end{algorithm}

%
%
\subsection{Complete pivoting}\label{ssec:fp}
%
%

To test the effectiveness of the described pivoting techniques, we implemented
complete pivoting in Algorithm~\ref{alg:GKO_TP}.
The overall computational cost is $O(n^3)$, since the algorithm requires at
step $k$ the reconstruction of a square submatrix of size $n-k+1$ by
formula~\eqref{eq:Cl}.
To keep storage linear with respect to $n$, we compute its columns one
at a time, saving the largest entry and its position in the column.

\begin{algorithm}[ht]
  \begin{algo}
    \REQUIRE $n,r,d \in\N^+$
    \REQUIRE $\bm{t},\bm{s}\in\C^n$; $G,H\in\C^{n\times r}$; $\bm{b}\in\C^{n\times d}$
    \STATE[] \phantom{\algofont{Require: }}s.t. $t_i\neq s_j$, $\forall i,j$, 
    and $s_i\neq s_j$, $\forall i\neq j$
    \STATE $H\leftarrow H^*$
    \STATE[3.9:]  $Q \leftarrow I_n$
    \FOR   $k=1$ \To $n$
    \FOR[z.1:]$j=k$ \To $n$
    \STATE[4.2:]  $\bm{v}_{k:n}\leftarrow(G_{k:n,:}\cdot H_{:,j})./(t_{k:n}-s_j)$
    \STATE[4.3:]  \textsc{find} $r_j,c_j$ \textsc{such that} $c_j:=|\bm{v}_{r_j}|=\max|\bm{v}_{k:n}|$
    \ENDFOR[4.4:]
    \STATE[4.5:]  \textsc{find} $i_2$ \textsc{such that} $c_{i_2}=\max(\bm{c}_{k:n})$
    \STATE[4.6:] $i_1\leftarrow r_{i_2}$
    \STATE[4.7:] \If $i_1\neq k$ \Then \textsc{swap rows} $k,i_1$ \textsc{of} $G$, $\bm{t}$, $\bm{b}$
    \STATE[4.8:] \If $i_2\neq k$ \Then \textsc{swap columns} $k,i_2$ \textsc{of}
    $H$, $\bm{s}$, $Q^T$
    \STATE[\dots]  \emph{\{ Insert lines 5--16 from Algorithm~\ref{alg:GKOstructured} \}}
    \algosetline{17}
    \ENDFOR
    \OUTPUT       $\bm{x}\leftarrow Q^T\bm{b}$
  \end{algo}
  \caption{Solver for $C\bm{x}=\bm{b}$ with complete pivoting.}
  \label{alg:GKO_TP}
\end{algorithm}

%
%
\subsection{Singularity detection}\label{ssec:sing}
%
%
A robust linear system solver should be able to detect singularity and to warn
the user about ill-conditioning.
Matlab \emph{backslash} operator, in the case of a square, non sparse matrix
$A$, performs the first task by checking if a diagonal element of the $U$
factor of $A$, computed by the \ttt{dgetrf} routine of LAPACK~\cite{lapack92},
is exactly zero.
Ill-conditioning is detected by estimating the 1-norm condition number of the
linear system by the \ttt{dgecon} routine of LAPACK, and a warning is issued if
the estimate \ttt{rcond} is less than $\ttt{eps}=2^{-52}$.

In our code we check the singularity in the same way, and we detect
ill-conditioning by computing the 1-norm condition number of the $U$ factor,
not to increase the overall computational load.
In fact, as factorization~\eqref{eq:fmb} shows, we compute the matrix $U$ by
rows (line 7 of Algorithm~\ref{alg:GKOstructured}) and $U^{-1}$ by columns
(line 5).
So, it is easy to obtain the 1-norm of $U^{-1}$, and to update step by step
the sums of the columns of $U$.
If the reciprocal of the resulting condition number is less than $\ttt{eps}$, a
warning is displayed.

%
%
\section{Package description}\label{sec:impl}
%
%
Our package is written for the most part in Matlab~\cite{matlab79}, and it has
been developed using version 7.9 (R2009b) on Linux, but we tested it on
various other versions, starting from 7.4.
Full documentation for every function in the package is accessible via the
Matlab \ttt{help} command, and the code itself is extensively commented.

The package is available at the web page 
\begin{verse}
\url{http://bugs.unica.it/~gppe/soft/}
\end{verse}
and it is distributed as a compressed archive file.
By uncompressing it, the directory \ttt{drsolve} which contains the software is
created.
This directory must be added to the Matlab search path, either by the
\ttt{addpath} command or using the menus available in the graphical user
interface.
Depending on the operating system and the Matlab version, some additional work
may be required; the details of the installation are discussed in the
\ttt{README.txt} file, located in the \ttt{drsolve} directory.

After the installation, the user should execute the script \ttt{validate.m} in
the subdirectory \ttt{drsolve/validate}.  It will check that the installation
is correct, that all the files work properly, and it will give some hints on
how to fix installation problems.

The core of the package is the function \ttt{clsolve}, which contains the
Matlab implementation of Algorithms~\ref{alg:GKOstructured}--\ref{alg:GKO_TP}.
Since these algorithms are heavily based on \emph{for} loops, which in general
degrade the performance of Matlab code, we reimplemented this function
in C language, using the Matlab C-MEX interface.
The MEX library allows to compile a Fortran or C subroutine, which can be
called by Matlab with the usual syntax.
To optimize the C code and to take advantage from the computer architecture, we
made an extensive use of the BLAS library~\cite{blas02}, keeping at a minimum
the use of explicit \emph{for} loops.
When the \ttt{clsolve} function is called, the compiled version is executed, if
it is available, otherwise the Matlab version is executed, issuing a warning
message.

The C code is contained in the \ttt{drsolve/src} subdirectory.
Its compilation is straightforward on Linux and Mac OS X (and we expect the
same behaviour on other Unix based platforms), while it is a bit more involved
on Matlab for Windows, which officially supports the MEX library only for some
commercial compilers.
Some further difficulties are due to the presence of complex variables in the
code, which are not supported by the minimal C compiler (\ttt{lcc}) distributed
with Matlab for Windows, and the need to link the BLAS and
LAPACK~\cite{lapack92} libraries.
We produced executables on Windows using a porting of the GNU-C
compiler~\cite{mingw} and the ``MEX configurer'' Gnumex~\cite{gnumex}.
Executables for various versions of Matlab on Linux, Mac OS X, and Windows, are
provided in our package; see the \texttt{README.txt} file for details on the
compilation process, and on the activation of the available executables.

An important technical remark which affects the 64-bit operating systems is the
following.
Starting from Matlab 7.8, the type of the integer variables used in the BLAS
and LAPACK libraries provided with Matlab changed from \ttt{int} to
\ttt{ptrdiff\_t}.
This has no effect on 32-bit architectures, but on 64-bit systems the size of
the variables changed from 4 to 8 bytes.
For this reason, to compile our C-MEX program \ttt{clsolve.c}, running Matlab
version 7.7 or less on a 64-bit system, it is necessary to uncomment one of the
first lines of the program, which is clearly highlighted in the code.

Besides the algorithms for solving Cauchy-like linear systems, coded in the
function \ttt{clsolve}, the package includes six simple interface programs,
listed in Table~\ref{tab:solver}, to solve linear systems with the displacement
structures discussed in Section~\ref{sec:displ}.
All these functions transform the system into a Cauchy like one by using the
corresponding conversion routine (see Table~\ref{tab:solver}), call
\ttt{clsolve} to compute the solution, and then recover the solution of the
initial system.

\begin{table}
\small
\centering
\begin{tabular}{lll}
\multicolumn{1}{c}{structure} & \multicolumn{1}{c}{solver} &
\multicolumn{1}{c}{conversion} \\
\hline
Cauchy-like & \ttt{clsolve} \\
Toeplitz & \ttt{tsolve} & \ttt{t2cl} \\
Toeplitz-like & \ttt{tlsolve} & \ttt{tl2cl} \\
Toeplitz+Hankel & \ttt{thsolve} & \ttt{th2cl} \\
Toeplitz+Hankel-like & \ttt{thlsolve} & \ttt{thl2cl} \\
Vandermonde & \ttt{vsolve} & \ttt{v2cl} \\
Vandermonde-like & \ttt{vlsolve} & \ttt{vl2cl} \\
\hline
\end{tabular}
\caption{Solvers and conversion routines.}
\label{tab:solver}
\end{table}
For example, given the Cauchy-like system $C\bm{x}=\bm{b}$, with displacement
$D_\bm{t}C-CD_\bm{s}=GH^*$, the commands
\begin{quote}
\begin{verbatim}
piv = 1;
x = clsolve(G,H,t,s,b,piv);
\end{verbatim}
\end{quote}
solve the system by Algorithm~\ref{alg:GKOssPP}, i.e., with partial pivoting.
The \ttt{piv} variable is used to select the solution algorithm.
To solve a Toeplitz linear system $T\bm{x}=\bm{b}$ with Gu's pivoting, one can
issue the following commands
\begin{quote}
\begin{verbatim}
piv = 4;
x = tsolve(c,r,b,piv);
\end{verbatim}
\end{quote}
where \ttt{c} and \ttt{r} denote the first column and the first row of $T$,
respectively.

The conversion routines and the test programs rely on some auxiliary routines,
listed in Table~\ref{tab:aux}.
The most relevant routines perform fast matrix products involving the unitary
matrices $F_\phi$, $\mathcal{S}$, and $\mathcal{C}$, introduced in
Section~\ref{sec:displ}.
The routines \ttt{stimes} and \ttt{ctimes}, concerning the Toeplitz+Hankel-like
structure, require the Matlab commands \ttt{dst} and \ttt{dct}/\ttt{idct},
available in the PDE Toolbox and in the Signal Processing Toolbox,
respectively.

\begin{table}[htb]
\small
\centering
\begin{tabular}{ll|ll}
\hline
\ttt{ftimes} & compute $F_\phi\bm{v}$ & \ttt{ttimes} & compute $T\bm{v}$ \\
\ttt{ctimes} & compute $\mathcal{C}\bm{v}$ & \ttt{cltimes} & compute $C\bm{v}$ \\
\ttt{stimes} & compute $\mathcal{S}\bm{v}$ & \ttt{cl2full} & assemble $C$ \\
\ttt{nroots1} & $n$th roots of $\e^{\ii\phi}$ \\
\hline
\end{tabular}
\caption{Auxiliary routines.}
\label{tab:aux}
\end{table}

The subdirectory \ttt{drsolve/test} contains the test programs which reproduce
the numerical experiments discussed in the next section.

%
%
\section{Numerical results}\label{sec:numres}

In this Section we present a selection of numerical experiments,
aimed to verify the effectiveness of our package, and to compare its
performance with other methods.

The numerical results were obtained with Matlab 7.9 (Linux 64-bit version), on
a single processor computer (AMD Athlon 64 3200+) with 1.5 Gbyte RAM, running
Debian GNU/Linux 5.0.
Each numerical experiment can be repeated by running the corresponding script
in the \ttt{drsolve/test} directory.
Except where explicitly noted, the C-MEX version of the \ttt{clsolve}
function was used.
In the experiments, the right hand side of each linear system corresponds
to the solution $(1,\ldots,1)^T$, errors are measured using the infinity norm,
and execution times are expressed in seconds.

\begin{table}[htb]
\small
\centering
\begin{tabular}{lcccc}
& disp.~rank & $\cond(A)$ & errors & exec.~time \\
\hline
Vandermonde & 1 & $5.0\cdot 10^{3}$ 
	& $4.3\cdot 10^{-13}/8.3\cdot 10^{-13\strut}$ & $0.80/26.52$ \\
Toeplitz & 2 & $3.5\cdot 10^{4}$ & $1.3\cdot 10^{-12}/4.6\cdot 10^{-12}$ 
	& $0.87/25.12$ \\
Toeplitz+Hankel & 4 & $5.0\cdot 10^{5}$ & $1.6\cdot 10^{-7}/1.5\cdot 10^{-11}$ 
	& $2.82/22.86$ \\
Cauchy-like & 5 & $4.5\cdot 10^{3}$ & $2.7\cdot 10^{-12}/5.5\cdot 10^{-13}$ 
	& $3.04/22.20$ \\
\hline
\end{tabular}
\caption{Errors and execution times in the solution of random complex 
linear systems of size 2048, belonging to four different structured classes.
The two figures in the rightmost columns are obtained by \ttt{drsolve} with
partial pivoting, and by Matlab \emph{backslash}, respectively.}
\label{tab:tab1}\end{table}

The first test concerns the solution of a complex linear system of size 2048,
whose matrix is either Vandermonde, Toeplitz, Toeplitz+Hankel, or Cauchy-like;
in the last case the displacement rank is 5.
The linear systems were constructed from random complex data; see the file
\ttt{test1.m} for details.
The errors and execution times, reported in Table~\ref{tab:tab1}, were obtained
by the solvers listed in~Table~\ref{tab:solver}, with partial pivoting, and by
Matlab \emph{backslash}.
It can be seen that in all cases the structured solver is much faster, as
expected, and that there is a significant loss in accuracy only in the
Toeplitz+Hankel case, where the condition number is larger.

\begin{figure}[htbp]
\centering\resizebox{.8\linewidth}{!}{\includegraphics{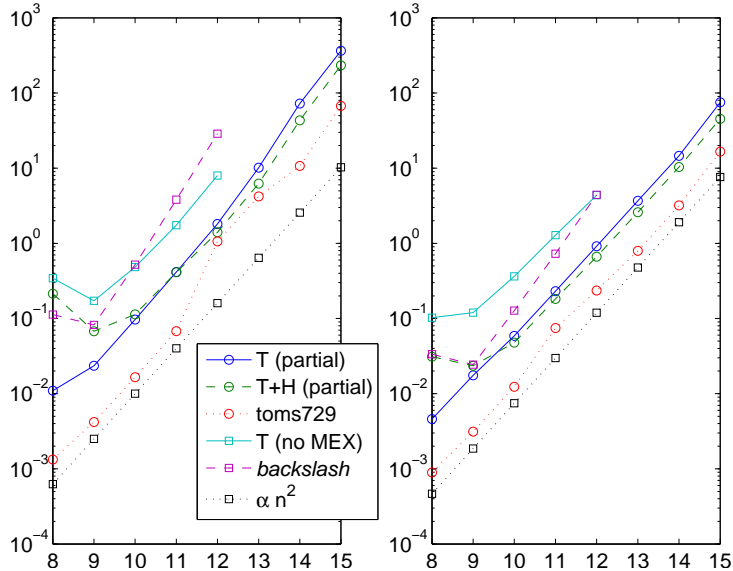}}
\caption{Execution time on a single processor computer (left) and on a
\emph{quad-core} computer (right) for a random real Toeplitz system of size
$2^k$, $k=8,\ldots,15$.}
\label{fig:time}
\end{figure}

In Figure~\ref{fig:time} the performance in terms of execution time is further
investigated.
In this case, a random real Toeplitz system of size $2^k$, $k=8,\ldots,15$, is
solved by \ttt{tsolve} and \ttt{thsolve} with partial pivoting (setting to zero
the Hankel part of the input matrix), and by the subroutine
\ttt{toms729}~\cite{hc92}, based on a look-ahead Levinson algorithm, with the
\ttt{pmax} parameter set to 10; a Matlab MEX gateway for this subroutine is
available~\cite{ar08}.
These methods are also compared, for $n\leq 4096$, to Matlab \emph{backslash},
and to a modified version of \ttt{tsolve} (the corresponding data are labelled
as ``no MEX''), which calls the Matlab implementation of \ttt{clsolve}.
The computation is repeated on a 4 processors computer (Intel Core2 Quad
Q6600), running the same operating system and Matlab version.
It results that \ttt{tsolve} is about 6 times faster when the compiled version
of \ttt{clsolve} is called.
Moreover, \ttt{toms729} is faster than our solvers, although it has the same
order of complexity.

\begin{figure}[htbp]
\centering\resizebox{.7\linewidth}{!}{\includegraphics{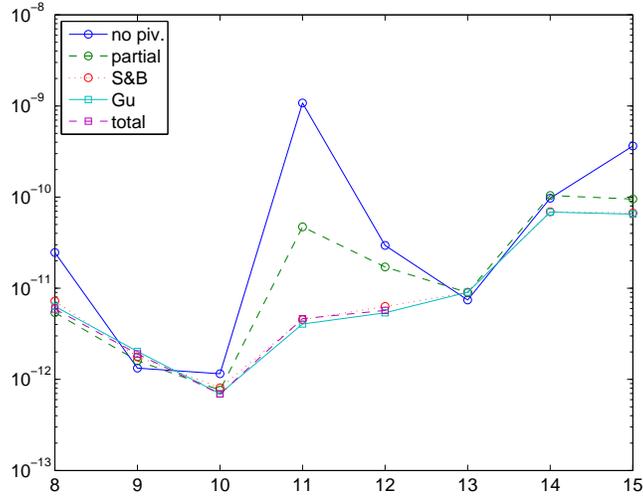}}
\caption{Error comparison between different pivoting techniques for a random
real Toeplitz system of size $2^k$, $k=8,\ldots,15$.}
\label{fig:erranpiv}
\end{figure}

Figure~\ref{fig:erranpiv} reports the errors in the solution of the same set of
real Toeplitz linear systems.  In this case, the different pivoting strategies
implemented in our package are compared; the results related to total pivoting
were computed only for $n\leq 4096$.
This figure shows that Sweet and Brent's pivoting (label ``S\&B'') and Gu's pivoting may produce
a very good approximation to total pivoting, leading to a significant
improvement in accuracy; we verified that the execution time is not
significantly affected by
the additional load required by these pivoting techniques.

\begin{figure}[htbp]
\centering\resizebox{.8\linewidth}{!}{\includegraphics{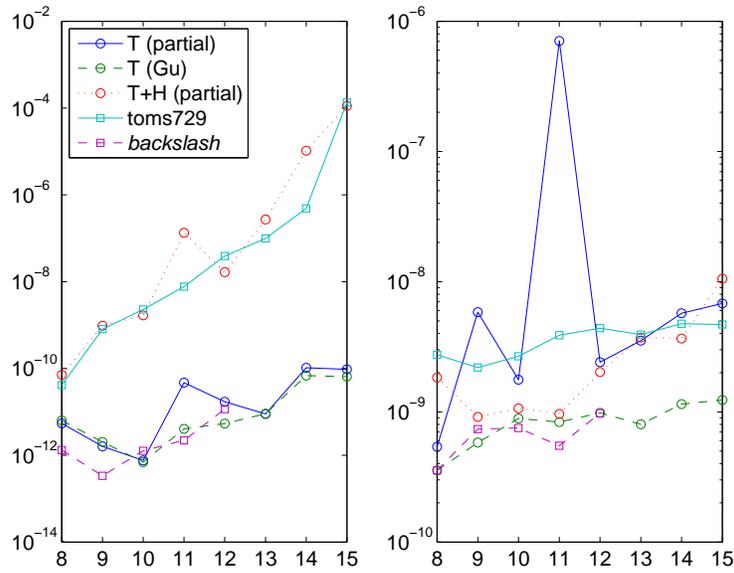}}
\caption{Error comparison between different solvers for a random real Toeplitz
system (left) and a Gaussian system (right) of size $2^k$, $k=8,\ldots,15$.}
\label{fig:ergauran}
\end{figure}

In Figure~\ref{fig:ergauran} we compare the accuracy of three of our Toeplitz
solvers, namely \ttt{tsolve} with partial and Gu's pivoting, and
\ttt{thsolve} with partial pivoting (setting to zero the Hankel part of the
matrix), to \ttt{toms729} and Matlab \emph{backslash}.
These methods are applied to the same set of Toeplitz linear systems used in
the previous experiments, and to Gaussian linear systems, whose matrix is
defined by
$$
a_{ij} = \sqrt{\frac{\sigma}{2\pi}} \e^{-\frac{\sigma}{2}(i-j)^2},
\quad \sigma = 0.3,
$$
and whose asymptotic condition number is $6.96\cdot 10^6$~\cite{ms05}.
The fact that Gu's pivoting produces very accurate results, generally
comparable to Matlab \emph{backslash}, is confirmed.
Moreover, the \ttt{toms729} and \ttt{thsolve} functions are often less
stable than the other methods, and there are cases in which partial pivoting
leads to a substantial error amplification.
The data displayed in Figures~\ref{fig:time}--\ref{fig:ergauran} were computed
by the scripts \ttt{test2.m} and \ttt{test3.m}.

\begin{figure}[htbp]
\centering\resizebox{.7\linewidth}{!}{\includegraphics{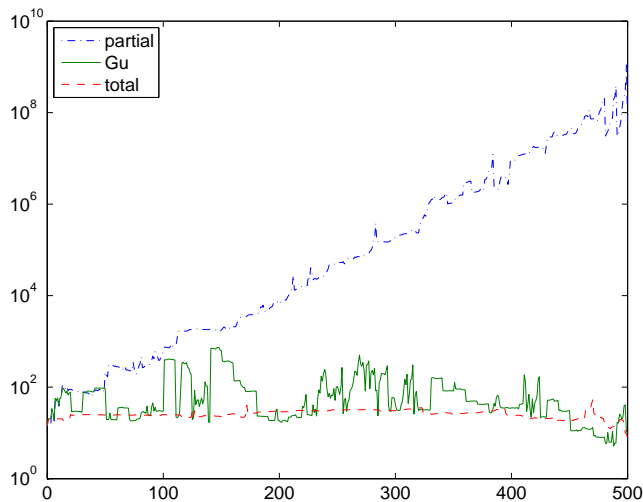}}
\caption{Growth of the right generator $H$ of the perturbed Bini-Boito example,
when three different pivoting techniques are applied.}
\label{fig:bbnormh}
\end{figure}

It is reported in the literature (see, e.g.,~\cite{sb95}) that the generalized
Schur algorithm may cause a large growth of the generators entries, causing the
error on the solution to be much larger than expected when standard Gaussian
elimination with partial pivoting is applied.
One such example is described in~\cite{bb10}, where the Sylvester matrix of two
polynomials is considered.
It is known that the rank deficiency of a Sylvester matrix equals the degree of
the polynomials GCD.
We consider a matrix of size 512, corresponding to two random polynomials with
a common factor of degree 20, and, since any Sylvester matrix is Toeplitz-like,
we convert it to a Cauchy-like matrix.
The generators of the matrix are then perturbed in order to make it
nonsingular, and the resulting system is solved by \ttt{clsolve}, with either
partial, Gu's or total pivoting; see \ttt{test4.m}.

The largest absolute value of the right generator entries at each iteration is
reported in Figure~\ref{fig:bbnormh}.
It is clear that, in this example, Gu's pivoting prevents generator growth as
much as total pivoting, while partial pivoting produces an exponential growth.
We remark that the norm of the solution errors corresponding to partial, Gu's,
total pivoting, and Matlab \emph{backslash}, are $1.1$, $1.1\cdot 10^{-5}$,
$2.5\cdot 10^{-6}$, and $3.5\cdot 10^{-6}$, respectively.

\begin{table}[htb]
\small
\centering
\begin{tabular}{lcccc}
& partial & Gu & total & \emph{backslash} \\
\hline
S\&B displacement & $3.3\cdot 10^{-3\strut}$ & $4.3\cdot 10^{-3}$ & $4.5\cdot 10^{-3}$ & $6.5\cdot 10^{-5}$ \\
alternative displ. & $1.0\cdot 10^{-14}$ & $1.3\cdot 10^{-14}$ & $4.6\cdot 10^{-14}$ & $7.7\cdot 10^{-15}$ \\
\hline
$\max_k|G_1^{(k)\strut}|/|G_1^{(0)}|$ & $ 1.6$ & $12.4$ & $ 1.6$ \\
$\max_k|H_1^{(k)}|/|H_1^{(0)}|$ & $11.3$ & $ 1.0$ & $11.3$ \\
\hline
\end{tabular}
\caption{Solution of the Sweet and Brent example, with two different choices of
generators.  
The first two lines reports the errors in the solution of the linear system by
\ttt{clsolve} (various pivoting) and \emph{backslash}.
For the first generators pair, the last two lines display the ratios between
the maximum entry of the generators computed at each iteration step, and the
maximum entry of the initial generators.}
\label{tab:tab5}\end{table}

We applied the same methods to the solution of an example proposed by Sweet and
Brent in~\cite{sb95}.
We consider a Cauchy-like matrix with knots on the unit circle, having the
following generators
\begin{equation}
G_1 = \begin{bmatrix} \bm{e} & \bm{e}+\tau\bm{f} \end{bmatrix}, \quad
H_1 = \begin{bmatrix} \bm{e} & -\bm{e} \end{bmatrix},
\label{eq:sbgen}
\end{equation}
with $\bm{e}=n^{-1/2}(1,\ldots,1)^T$, $\bm{f}=n^{-1/2}(-1,1,\ldots,(-1)^n)^T$,
and $\tau=10^{-12}$.
These generators produce huge cancellations when the elements of the matrix are
recovered.
To better investigate this example, we also considered the generators
$$
G_2 = -\tau\bm{f}, \quad
H_2 = \bm{e},
$$
which give an equivalent representation of the matrix, preventing
cancellation; see \ttt{test5.m}.
From Table~\ref{tab:tab5}, it results that there is no growth
associated to the generators~\eqref{eq:sbgen} given by Sweet and Brent, and
that the loss of accuracy is mainly due to the cancellation in the product
$G_1H_1^*$.

\begin{figure}[htbp]
\centering\resizebox{.7\linewidth}{!}{\includegraphics{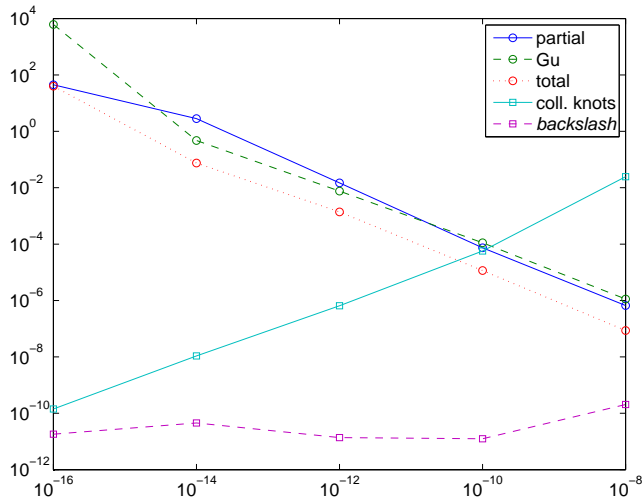}}
\caption{Error in the solution of a Cauchy-like linear system of size 260 with
\emph{almost} multiple knots.  Each knot is repeated 5 times, with a different
random perturbation producing a relative error less than~$\tau$; the values of
$\tau$ are reported on the horizontal axis.}
\label{fig:clmk}
\end{figure}

Algorithm~\ref{alg:GKOss} can be applied to a Cauchy-like matrix with
multiple knots; see Section~\ref{ssec:ss}.
We verified that its numerical performance is not influenced by the knots
multiplicity and that, when $s_i\neq s_j$, $\forall i\neq j$, it is roughly as
accurate as Algorithm~\ref{alg:GKOstructured}, so we prefer
to investigate here the particular situation of \emph{almost} multiple knots.
In the script \ttt{test6.m}, we construct a Cauchy-like matrix of size $260$
and displacement rank 5, with random generators.
Its knots are obtained starting from 52 equispaced points on the unit circle,
and replicating them 4 times, each time with a perturbation producing a
relative error less than~$\tau$.
Figure~\ref{fig:clmk} reports the errors obtained solving the resulting
linear system for $\tau=10^{-8},10^{-10},\ldots,10^{-16}$.
Our \ttt{clsolve} function was applied with various pivoting techniques; the
results labelled as ``coll.~knots'' were obtained by collapsing the knots,
i.e., removing the perturbation $\tau$, and then solving the system using
Algorithm~\ref{alg:GKOssPP}, which includes partial pivoting.
In this way, a different system is solved, since the matrix is perturbed, but
when $\tau$ tends to zero this approach is much more accurate than
Algorithm~\ref{alg:GKOstructured}.
The condition number of the matrix is about $10^5$ for any $\tau$; this fact
is confirmed by the results obtained by Matlab \emph{backslash}.

\begin{figure}[htbp]
\centering\resizebox{.7\linewidth}{!}{\includegraphics{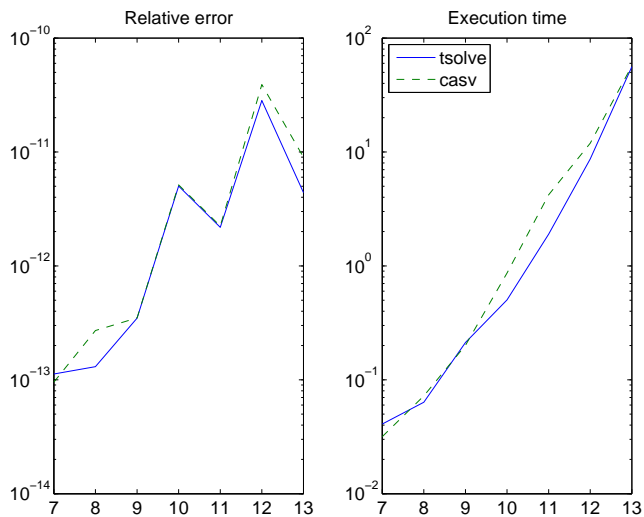}}
\caption{Comparison between \ttt{tsolve} with partial pivoting and the
function \ttt{casv} from~\protect\cite{pol10}. Both methods are applied to
random real Toeplitz systems of size $2^k$, $k=7,\ldots,13$.}
\label{fig:tsolcasv}
\end{figure}

While working on this paper, we became aware of another approach to apply the
generalized Schur algorithm to a Cauchy-like matrix, using partial pivoting,
and requiring $O(n)$ memory locations~\cite{pol10}.
As the author kindly sent us the code developed for his paper, we report a
comparison of his method and \ttt{tsolve} with partial pivoting.
Figure~\ref{fig:tsolcasv} reports the errors and the execution times obtained
by applying the Matlab version of both methods to the solution of random real
Toeplitz linear systems; the results were computed by suitably modifying
the script \ttt{test2.m}.
It results that the algorithms are essentially equivalent.
Performing the same test using the compiled versions of both methods produces
the same errors, but our implementation appears to run faster.
We believe that this is due to a different level of code optimization.

\section*{Acknowledgment}

We thank Mario Arioli, Marc Van Barel, and Federico Poloni, for fruitful
discussions and suggestions.
We also thank Paola Boito for providing us the software to produce one of the
numerical examples, and Federico Poloni for sending us the function \ttt{casv}
from~\cite{pol10}.


\begin{thebibliography}{10}

\bibitem{lapack92}
E.~Anderson, Z.~Bai, C.~Bischof, J.~Demmel, J.~Dongarra, J.~Du~Croz,
  A.~Greenbaum, S.~Hammarling, A.~McKenney, S.~Ostrouchov, and D.~Sorensen.
\newblock {\em LAPACK Users' Guide}.
\newblock SIAM, Philadelphia, 1992.

\bibitem{ar08}
A.~Aric\`o and G.~Rodriguez.
\newblock {\em \texttt{toms729gw}: a Matlab (Fortran) MEX Gateway for TOMS
  Algorithm 729, by P. C. Hansen and T. Chan}.
\newblock University of Cagliari, 2008.
\newblock \newline Available at: \url{http://bugs.unica.it/~gppe/soft/}.

\bibitem{toepack83}
O.~B. Arushanian, M.~K. Samarin, V.~V. Voevodin, E.~Tyrtyshnikov, B.~S. Garbow,
  J.~M. Boyle, W.~R. Cowell, and K.~W. Dritz.
\newblock The {TOEPLITZ} {P}ackage {U}sers' {G}uide.
\newblock Technical Report ANL-83-16, Argonne National Laboratory, 1983.

\bibitem{bgs70}
R.~H. Bartels, G.~H. Golub, and M.~A. Saunders.
\newblock Numerical techniques in mathematical programming.
\newblock In {\em Nonlinear {P}rogramming ({P}roc. {S}ympos., {U}niv. of
  {W}isconsin, {M}adison, {W}is., 1970)}, pages 123--176. Academic Press, New
  York, 1970.

\bibitem{bb10}
D.A. Bini and P.~Boito.
\newblock A fast algorithm for approximate polynomial gcd based on structured
  matrix computations.
\newblock In {\em Numerical Methods for Structured Matrices and Applications:
  Georg Heinig memorial volume}. Birkh\"auser, Basel, 2010.
\newblock To appear.

\bibitem{bjo92}
{\AA}.~Bj{\"o}rck.
\newblock Pivoting and stability in the augmented system method.
\newblock In {\em Numerical Analysis 1991 ({D}undee, 1991)}, volume 260 of {\em
  Pitman Res. Notes Math. Ser.}, pages 1--16. Longman Sci. Tech., Harlow, 1992.

\bibitem{blas02}
L.~S. Blackford, J.~Demmel, J.~Dongarra, I.~Duff, S.~Hammarling, G.~Henry,
  M.~Heroux, and L.~Kaufman.
\newblock An updated set of basic linear algebra subprograms ({BLAS}).
\newblock {\em ACM Trans. Math. Software}, 28(2):135--151, 2002.

\bibitem{cr97}
D.~Calvetti and L.~Reichel.
\newblock Factorizations of {C}auchy matrices.
\newblock {\em J. Comput. Appl. Math.}, 86(1):103--123, 1997.

\bibitem{fmkl79}
B.~Friedlander, M.~Morf, T.~Kailath, and L.~Ljung.
\newblock New inversion formulas for matrices classified in terms of their
  distance from {T}oeplitz matrices.
\newblock {\em Linear Algebra Appl.}, 27:31--60, 1979.

\bibitem{gko95}
I.~Gohberg, T.~Kailath, and V.~Olshevsky.
\newblock Fast {G}aussian elimination with partial pivoting for matrices with
  displacement structure.
\newblock {\em Math. Comp.}, 64(212):1557--1576, 1995.

\bibitem{gu98}
M.~Gu.
\newblock Stable and efficient algorithms for structured systems of linear
  equations.
\newblock {\em SIAM J. Matrix Anal. Appl.}, 19(2):279--306, 1998.

\bibitem{hc92}
P.~C. Hansen and T.~F. Chan.
\newblock Fortran subroutines for general {T}oeplitz systems.
\newblock {\em ACM Trans. Math. Software}, 18(3):256--273, 1992.

\bibitem{hei95a}
G.~Heinig.
\newblock Inversion of generalized {C}auchy matrices and other classes of
  structured matrices.
\newblock In {\em Linear Algebra for Signal Processing ({M}inneapolis, {MN},
  1992)}, volume~69 of {\em IMA Vol. Math. Appl.}, pages 63--81. Springer, New
  York, 1995.

\bibitem{hb97}
G.~Heinig and A.~Bojanczyk.
\newblock Transformation techniques for {T}oeplitz and {T}oeplitz-plus-{H}ankel
  matrices. {I}. {T}ransformations.
\newblock {\em Linear Algebra Appl.}, 254:193--226, 1997.

\bibitem{hb98}
G.~Heinig and A.~Bojanczyk.
\newblock Transformation techniques for {T}oeplitz and {T}oeplitz-plus-{H}ankel
  matrices. {I}{I}. {A}lgorithms.
\newblock {\em Linear Algebra Appl.}, 278(1-3):11--36, 1998.

\bibitem{hr84}
G.~Heinig and K.~Rost.
\newblock {\em Algebraic Methods for {T}oeplitz-Like Matrices and Operators},
  volume~13 of {\em Operator Theory: Advances and Applications}.
\newblock Birkh\"auser Verlag, Basel, 1984.

\bibitem{kc94}
T.~Kailath and J.~Chun.
\newblock Generalized displacement structure for block-{T}oeplitz,
  {T}oeplitz-block, and {T}oeplitz-derived matrices.
\newblock {\em SIAM J. Matrix Anal. Appl.}, 15(1):114--128, 1994.

\bibitem{kkm79b}
T.~Kailath, S.~Y. Kung, and M.~Morf.
\newblock Displacement ranks of matrices and linear equations.
\newblock {\em J. Math. Anal. Appl.}, 68(2):395--407, 1979.

\bibitem{ks95}
T.~Kailath and A.~H. Sayed.
\newblock Displacement structure: theory and applications.
\newblock {\em SIAM Rev.}, 37(3):297--386, 1995.

\bibitem{ks99}
T.~Kailath and A.~H. Sayed, editors.
\newblock {\em Fast Reliable Algorithms for Matrices with Structure},
  Philadelphia, PA, 1999. Society for Industrial and Applied Mathematics
  (SIAM).

\bibitem{maseteam}
Katholieke Universiteit Leuven, Department of Computer Science.
\newblock {\em MaSe-Team (Matrices having Structure)}, 2010.
\newblock \newline Available at:
  \url{http://www.cs.kuleuven.ac.be/~marc/software/}.

\bibitem{matlab79}
The MathWorks, Natick.
\newblock {\em Matlab ver.\ 7.9}, 2009.

\bibitem{ms05}
C.V.M. van~der Mee and S.~Seatzu.
\newblock A method for generating infinite positive self-adjoint test matrices
  and {R}iesz bases.
\newblock {\em SIAM Journal on Matrix Analysis and Applications},
  26(4):1132--1149, 2005.

\bibitem{pol10}
F.~Poloni.
\newblock A note on the {$O(n)$}-storage implementation of the {GKO} algorithm
  and its adaptation to {T}rummer-like matrices.
\newblock {\em Numer. Algorithms}, 2010.
\newblock To appear.

\bibitem{rod06}
G.~Rodriguez.
\newblock Fast solution of {T}oeplitz- and {C}auchy-like least-squares
  problems.
\newblock {\em SIAM J. Matrix Anal. Appl.}, 28(3):724--748 (electronic), 2006.

\bibitem{sie65}
I.~H. Siegel.
\newblock Deferment of computation in the method of least squares.
\newblock {\em Math. Comp.}, 19:329--331, 1965.

\bibitem{gnumex}
SourceForge.net.
\newblock {\em Gnumex}, 2000.
\newblock \newline Available at: \url{http://gnumex.sourceforge.net/}.

\bibitem{mingw}
SourceForge.net.
\newblock {\em MinGW}, 2008.
\newblock \newline Available at: \url{http://www.mingw.org/}.

\bibitem{sb95}
D.~R. Sweet and R.~P. Brent.
\newblock Error analysis of a fast partial pivoting method for structured
  matrices.
\newblock In F.~T. Luk, editor, {\em Advanced Signal Processing Algorithms},
  volume 2563, pages 266--280, San Diego, 1995. SPIE.

\bibitem{bezout}
University of Pisa, Department of Mathematics.
\newblock {\em Structured matrix analysis: numerical methods and applications},
  2010.
\newblock \newline Available at: \url{http://bezout.dm.unipi.it/}.

\end{thebibliography}

\end{document}